\documentclass[12pt]{article}
\usepackage{amsmath}
\usepackage{amssymb}
\usepackage{vector}
\def\eps{\varepsilon}
\font\tencmmib=cmmib10 \skewchar\tencmmib '60
\newfam\cmmibfam
\textfont\cmmibfam=\tencmmib

\def\bbox{\quad\hbox{\vrule \vbox{\hrule \vskip2pt \hbox{\hskip2pt
\vbox{\hsize=1pt}\hskip2pt} \vskip2pt\hrule}\vrule}}
\def\lessim{\ \lower4pt\hbox{$
\buildrel{\displaystyle <}\over\sim$}\ }
\def\gessim{\ \lower4pt\hbox{$\buildrel{\displaystyle >}
\over\sim$}\ }

\def\E{{\cal E}}

\def\R{{\cal R}}

\def\eps{{\varepsilon}}

\def\la{{\Bigl\langle}}
\def\ra{{\Bigr\rangle}}
\def\hs{{\hat{\sigma}}}
\def\hS{{\hat{\vsi}}}
\def\hR{{\hat{R}}}
\def\pepsi{{p_{\eps}}}
\def\bq{{\bar{q}}}
\def\bra{\bigr\rangle^{-}}
\def\bla{{\bigl\langle}}

\def\qed{\hfill\break\rightline{$\bbox$}}
\newcommand{\e}{\mathbb{E}}
\newcommand{\p}{\mathbb{P}}
\newcommand{\Reals}{\mathbb{R}}

\newcommand{\vsi}{{\vec{\sigma}}}

\newtheorem{lemma}{Lemma}
\newtheorem{theorem}{Theorem}

\makeatletter
\@addtoreset{equation}{section}

\makeatother

\font\tencmmib=cmmib10 \skewchar\tencmmib '60
\newfam\cmmibfam
\textfont\cmmibfam=\tencmmib

\def\bbox{\quad\hbox{\vrule \vbox{\hrule \vskip2pt \hbox{\hskip2pt
\vbox{\hsize=1pt}\hskip2pt} \vskip2pt\hrule}\vrule}}
\def\lessim{\ \lower4pt\hbox{$
\buildrel{\displaystyle <}\over\sim$}\ }
\def\gessim{\ \lower4pt\hbox{$\buildrel{\displaystyle >}
\over\sim$}\ }

\def\eps{\varepsilon}

\def\go0{\to 0}

\def\la{\langle}

\def\leftitem#1{\item{\hbox to\parindent{\enspace#1\hfill}}}

\def\qed{\hfill\break\rightline{$\bbox$}}

\def\ra{\rangle}

\def\sg{\sigma}

\def\sg2{\sigma^2}

\def\__{_{\infty}}

\begin{document}

\title{
Cavity method in the spherical Sherrington-Kirkpatrick model.}

\author{ 
Dmitry Panchenko\thanks{Department of Mathematics, Massachusetts Institute
of Technology, 77 Massachusetts Ave, Cambridge, MA 02139
email: panchenk@math.mit.edu. Partially supported by NSF grant.}
}
\date{}

\maketitle

\begin{abstract}
We develop a cavity method in the spherical Sherrington-Kirkpatrick
model at high temperature and small external field. As one application
we compute the limit of the covariance matrix for fluctuations of the
overlap and magnetization.
\end{abstract}
\vspace{0.5cm}

Key words: Sherrington-Kirkpatrick model, cavity method.

Mathematics Subject Classification: 60K35, 82B44

\section{Introduction.}

The cavity method in the Sherrington-Kirkpatrick model \cite{SK} 
as described, for example, in Chapter 2 of \cite{SG}, is one of the most
important tools used to analyze the model in the high temperature region.
As a typical applications of the cavity method one can show that
the overlap of two spin configurations is nearly constant and its fluctuations
are Gaussian (see \cite{TAOP} or \cite{GT}).
When we tried to understand how the cavity method would look like in
the spherical SK model, the task turned out to be much
more difficult than expected, mostly, due to the fact that uniform measure on the sphere
is not a product measure. (Of course, we can not even compare this difficulty
with the real difficulty of discovering original cavity method in the classical SK
model.) As an applications, we study fluctuations of the overlap and magnetization
and compute their covariance matrix in the thermodynamic limit. We stop short of
proving a central limit theorem since our goal is to provide a reasonably simple illustration
of the cavity method. 

We consider a spherical SK model with Gaussian Hamiltonian
$H_N(\vsi)$ indexed by {\it spin configurations} $\vsi$ on the sphere
$S_N$ of radius $\sqrt{N}$ in $\Reals^N$. We will assume that
\begin{equation}
\frac{1}{N}\e H_N(\vsi^1)H_N(\vsi^2)=\xi(R_{1,2})
\label{cov}
\end{equation}
where $R_{1,2}=N^{-1}\sum_{i\leq N} \sigma_i^1\sigma_i^2$ is {\it the overlap} of 
configurations $\vsi^1, \vsi^2\in S_N$ and where the function
$\xi(x)$ is three times continuously differentiable. 
This model was studied in \cite{CS} and rigorously in \cite{sph}. Under the additional
assumptions on $\xi,$ 
\begin{equation}
\xi(0)=0,\, \xi(x)=\xi(-x),\,\xi''(x)>0 \mbox{ if } x>0,
\label{xi}
\end{equation}
the limit of the free energy
\begin{equation}
F_N=\frac{1}{N}\e \log \int_{S_N}\exp\Bigl(
\beta H_N(\vsi)+h\sum_{i\leq N}\sigma_i
\Bigr) \lambda_N(\vsi)
\label{FE}
\end{equation}
was computed in \cite{sph} for arbitrary {\it inverse temperature} $\beta>0$ and 
{\it external field} $h\in\Reals.$ Here $\lambda_N$ denotes the uniform 
probability measure on $S_N.$  
 
The main results of the present paper will be proved for small enough parameters
$\beta$ and $h,$ i.e. for very high temperature and small external field,
and {\it without} the assumptions in (\ref{xi}), i.e. not only for even spin
interactions. However, to motivate these
results we will first describe some implications of the results in \cite{sph}
that were proved under (\ref{xi}).
 
For small $\beta$ and $h$ the results in \cite{sph} imply that under (\ref{xi})
the limit of the free energy takes a particularly simple form:
\begin{equation}
\lim_{N\to\infty} F_N=\inf_{q\in[0,1]}\frac{1}{2}\Bigl(
h^2(1-q)+\frac{q}{1-q}+\log(1-q)+\beta^2\xi(1)-\beta^2\xi(q)
\Bigr).
\label{Flim}
\end{equation}
In fact, the entire {\it replica symmetric} region of parameters $\beta, h$
where (\ref{Flim}) holds
can be easily described using Proposition 2.1 in \cite{sph}.
The critical point equation for the infimum on the right hand side of (\ref{Flim}) is
\begin{equation}
h^2+\beta^2\xi'(q)=\frac{q}{(1-q)^2}.
\label{critq}
\end{equation}
For small enough $\beta$ the infimum in (\ref{Flim}) is achieved at $q=0$ if $h=0$
and at the unique solution $q$ of (\ref{critq}) if $h\not =0.$ 
Theorem 1.2 in \cite{sph} suggests that the distribution of the overlap 
$R_{1,2}$ with respect to the Gibbs measure is concentrated near $q$
and by analogy with the Ising SK model (see Chapter 2 in \cite{SG} or \cite{GT}) 
one expects that the distribution of $\sqrt{N}(R_{1,2} -q)$ is approximately Gaussian. 
The proof of this result in \cite{SG} was based on the {\it cavity method} and 
the main goal of the present paper is to develop the analogue of the cavity
method for the spherical SK model. As we shall see, the cavity method
for the spherical model will be much more involved due to the fact that the
measure $\lambda_N$ on the sphere $S_N$ is not a product measure and
it will take some effort to decouple one coordinate from the others. 
The "cavity computations" will also be more involved and instead of
proving a central limit theorem for the overlap we will only carry out the 
computation of the variance of $\sqrt{N}(R_{1,2}-q)$ and other related quantities.
Without doubt, with extra work the cavity method developed in this paper would yield
central limit theorems as well. 

It is interesting to note that our results imply the analogue of (\ref{Flim})
without the assumption (\ref{xi}). Namely, since we will prove that
for small $\beta$ and $h$ the overlap $R_{1,2}$ is concentrated near 
the unique solution $q$ of (\ref{critq}), it is a simple exercise to show
that in this case
\begin{equation}
\lim_{N\to\infty} F_N=\frac{1}{2}\Bigl(
h^2(1-q)+\frac{q}{1-q}+\log(1-q)+\beta^2\xi(1)-\beta^2\xi(q)
\Bigr).
\label{Flim2}
\end{equation}
To prove this, one only needs to compare the derivatives of both
sides with respect to $\beta$ since
$$
\frac{\partial F_N}{\partial \beta}=\beta(\xi(1)-\e\la\xi(R_{1,2})\ra).
$$
The rest of the paper is organized as follows. Main object of the
paper - the cavity interpolation - is presented in the next section 
where we also state its main properties such as control of the derivative 
and a way to compute certain moments at the end of the interpolation.
In Section 3 we describe our main application of the cavity method
- the so called second moment computations, which constitute the first
step toward proving the central limit theorems for the overlap and magnetization.
Most of the technical proofs are left until Sections 5 and 6.

\section{Cavity method.}

For certainty, from now on we assume that $h\not= 0$ and 
$\beta$ is smal enough so that $q$ is the {\it unique} solution of
(\ref{critq}). All the results below are proved without the assumption (\ref{xi}).
Given a configuration $\vsi\in S_N,$ we will denote $\eps=\sigma_N$
and for $i\leq N-1$ denote 
$$
\hs_i=\sigma_i\Bigr/\sqrt{\frac{N-\eps^2}{N-1}},
$$
so that a vector $\hS=(\hs_1,\ldots,\hs_{N-1})\in S_{N-1},$ i.e. 
$|\hS|=\sqrt{N-1}.$ We consider a Gaussian Hamiltonian $H_{N-1}(\hS)$ 
independent of $H_N(\vsi)$ such that
\begin{equation}
\frac{1}{N-1}\e H_{N-1}(\hS^1) H_{N-1}(\hS^2)=\xi(\hR_{1,2}),
\label{covhat}
\end{equation}
where $\hR_{1,2}=(N-1)^{-1}\sum_{i\leq N-1} \hs_i^1 \hs_i^2.$
We define an interpolating Hamiltonian by
\begin{eqnarray}
H_t(\vsi)
&=&
\sqrt{t}\beta H_N(\vsi)+\sqrt{1-t}\beta H_{N-1}(\hS) + h\sum_{i\leq N-1} \hs_i
\Bigl(1+t\Bigl(\sqrt{\frac{N-\eps^2}{N-1}}-1\Bigr)\Bigr) 
\nonumber
\\
&&
+ h\eps
+\sqrt{1-t}\eps z\beta\sqrt{\xi'(q)}-\frac{1}{2}(1-t)\eps^2 b
\label{Hamint}
\end{eqnarray}
where $z$ is a Gaussian r.v. independent of $H_N$ and $H_{N-1}$ and 
\begin{equation} 
b=h^2(1-q)+\beta^2(1-q)\xi'(q).
\label{b}
\end{equation}
The main idea in this interpolation (which was hardest to discover)
is that we interpolate directly between spin configurations on
$S_{N}$ and $S_{N-1}$!
The cavity Hamiltonian at $t=0$ is
\begin{equation}
H_0(\vsi)=\beta H_{N-1}(\hS)+h\sum_{i\leq N-1}\hs_i + \eps a -\frac{1}{2}\eps^2 b,
\label{Hamzero}
\end{equation}
where we introduced the notation
\begin{equation}
a=z\beta\sqrt{\xi'(q)}+h.
\label{a}
\end{equation}
The terms that do not depend on $\eps$ depend on the rest of the coordinates only through
$\hS\in S_{N-1}$ and, therefore, the Gibbs' average at $t=0$ for functions of the type
$f_1(\hS)f_2(\eps)$ will decouple, which is a crucial feature of the cavity method.
Another feature that one expects from this interpolation is that, as we will show,
along the interpolation annealed Gibbs averages do not change much. To show this, we will
first compute the derivatives along the interpolation.
Define 
$$
Z_t=\int_{S_N} \exp H_t(\vsi) d\lambda_N(\vsi)
$$ 
and for a function
$f:S_N^n\to\Reals$ define the Gibbs average of $f$ with respect to the 
Hamiltonian (\ref{Hamint}) by
\begin{equation}
\la f\ra_t=\frac{1}{Z_t^n}\int_{S_N^n}\exp\sum_{l\leq n} H_t(\vsi^l) d\lambda_N^n.
\label{Gibbst}
\end{equation}
Let $\nu_t(t)=\e \la f\ra_t.$ For $q$ in (\ref{critq}) we define
\begin{equation}
r=h(1-q).
\label{r}
\end{equation}
Let $\hR = (N-1)^{-1}\sum_{i\leq N-1}\hs_i.$ We define $a_l$ and $a_{l,l'}$ by
\begin{equation}
a_l=1-\eps_l^2,\,\,\,\,
2a_{l,l'}=\xi'(q)-\frac{1}{2}(\eps_l^2+\eps_{l'}^2)\left(q\xi''(q)+\xi'(q)\right)
+\eps_l\eps_{l'}\xi''(q).
\label{ais}
\end{equation}
The following holds.

\begin{theorem}\label{Thder}
We have
\begin{eqnarray}
\nu_t'(f)
&=&
\frac{h}{2}\sum_{l\leq n}\nu_t(f a_l(\hR_l - r))
-n\frac{h}{2}\nu_t(f a_{n+1} (\hR_{n+1}-r))
\nonumber
\\
&&+\
2\beta^2 \sum_{1\leq l<l'\leq n} \nu_t(f a_{l,l'}(\hR_{l,l'}-q))
-2n\beta^2 \sum_{l\leq n} \nu_t(f a_{l,n+1}(\hR_{l,n+1}-q))
\nonumber
\\
&&+\
n(n+1)\beta^2 \nu_t(f a_{n+1,n+2}(\hR_{n+1,n+2}-q))
+\nu_t(f\R),
\label{der}
\end{eqnarray}
where the remainder $\R$ is bounded by
$$
|\R|\leq \frac{L}{N}(\beta^2+h)\Bigl(1+\sum_{l\leq n+2}\eps_l^4\Bigr) 
+ L\beta^2 \sum_{1\leq l\not=l'\leq n+2}(1+\eps_l^2)(\hR_{l,l'}-q)^2.
$$
\end{theorem}
\textbf{Proof.} 
We start by writing
\begin{equation}
\nu_t'(f)=\e\Bigl\la f\sum_{l\leq n} \frac{\partial}{\partial t} H_t(\vsi^l)\Bigr\ra_t
-n \e\Bigl\la f\frac{\partial}{\partial t} H_t(\vsi^{n+1})\Bigr\ra_t
\label{der1}
\end{equation}
and
\begin{eqnarray}
\frac{\partial}{\partial t} H_t(\vsi)
&=&
\frac{\beta}{2\sqrt{t}}H_N(\vsi)-\frac{\beta}{2\sqrt{1-t}}H_{N-1}(\hS) + h\sum_{i\leq N-1} \hs_i
\Bigl(\sqrt{\frac{N-\eps^2}{N-1}}-1\Bigr) 
\nonumber
\\
&&
-\frac{1}{2\sqrt{1-t}}\eps z\beta\sqrt{\xi'(q)}+\frac{1}{2}\eps^2 b.
\label{Hamintder}
\end{eqnarray}
In order to use a Gaussian integration by parts (see, for example, (A.41) in \cite{SG}) 
we first compute the covariance
$$
\mbox{Cov}\Bigl(
H_t(\vsi^1),\frac{\partial}{\partial t} H_t(\vsi^{2})
\Bigr)
=\frac{\beta^2}{2}\Bigl(N\xi(R_{1,2})-(N-1)\xi(\hR_{1,2})
-\eps_1 \eps_2 \xi'(q)\Bigr),
$$
by (\ref{cov}) and (\ref{covhat}). We will rewrite this using Taylor's expansion
of $\xi(R_{1,2})$ near $\hR_{1,2}.$ We will use that
\begin{equation}
R_{1,2}=\hR_{1,2}+s(\eps_1,\eps_2)\hR_{1,2}
+N^{-1}\eps_1\eps_2
\label{ars}
\end{equation}
where
$$
s(\eps_1,\eps_2)=
\sqrt{\Bigl(1-\frac{\eps_1^2}{N}\Bigr)\Bigl(1-\frac{\eps_2^2}{N}\Bigr)} -1.
$$
Since
\begin{equation}
\Bigl|\sqrt{1+x}-1-\frac{x}{2}\Bigr|\leq L x^2 \mbox{ for } x\in[-1,1]
\label{xdel}
\end{equation}
we have
\begin{equation}
\Bigl|
s(\eps_1,\eps_2)+\frac{1}{2N}(\eps_1^2 +\eps_2^2)
\Bigr|
\leq \frac{L}{N^2}(\eps_1^4+\eps_2^4).
\label{epsils}
\end{equation}
By assumption, $\xi$ is three times continuously differentiable and
(\ref{ars}), (\ref{epsils}) imply 
$$
\Bigl|
\xi(R_{1,2})-\xi(\hR_{1,2})
-\xi'(\hR_{1,2})(R_{1,2}-\hR_{1,2})
\Bigl|\leq
\frac{L}{N^2}(\eps_1^4+\eps_2^4)
$$
and
$$
\Bigl|
\xi(R_{1,2})-\xi(\hR_{1,2})
+\frac{1}{2N}(\eps_1^2+\eps_2^2)\hR_{1,2}\xi'(\hR_{1,2})
-\frac{1}{N}\eps_1\eps_2\xi'(\hR_{1,2})
\Bigl|\leq
\frac{L}{N^2}(\eps_1^4+\eps_2^4).
$$
Therefore,
\begin{equation}
N \xi(R_{1,2}) -(N-1)\xi(\hR_{1,2})= \xi(\hR_{1,2}) 
-\frac{1}{2}(\eps_1^2+\eps_2^2)\hR_{1,2}\xi'(\hR_{1,2})
+ \eps_1\eps_2 \xi'(\hR_{1,2})+\R_1
\label{weget1}
\end{equation}
where from now on $\R_1$ will denote a quantity such that
$$
|\R_1| \leq \frac{L}{N}\Bigl(1+\sum_{l\leq n+2}\eps_l^4\Bigr).
$$ 
Since $\xi$ is three times continuously differentiable,
\begin{eqnarray*}
&&
\xi(\hR_{1,2})-\xi(q)=\xi'(q)(\hR_{1,2}-q)+\R_2,\,\,\,
\xi'(\hR_{1,2})-\xi'(q)=\xi''(q)(\hR_{1,2}-q)+\R_2,\,\,\,
\\
&&
\hR_{1,2}\xi'(\hR_{1,2})-q\xi'(q)=(\xi'(q)+q\xi''(q))(\hR_{1,2}-q)+\R_2,\,\,\,
\end{eqnarray*}
where $\R_2$ denotes a quantity such that 
$$
|\R_2|\leq L(\hR_{1,2}-q)^2.
$$
Using this in (\ref{weget1}) and recalling the definition of $a_{l,l'}$
in (\ref{ais}) we get
\begin{equation}
\mbox{Cov}\Bigl(
H_t(\vsi^l),\frac{\partial}{\partial t} H_t(\vsi^{l'})
\Bigr)
=\frac{\beta^2}{2}\Bigl(
2a_{l,l'}(\hR_{l,l'}-q)-\frac{1}{2}(\eps_l^2+\eps_{l'}^2)q\xi'(q)
+\xi(q) \Bigr) +\beta^2 \R_3
\label{cov1}
\end{equation}
where 
$$
|\R_3|\leq \frac{L}{N}\Bigl(1+\sum_{l\leq n+2}\eps_l^4\Bigr) 
+ L\sum_{l\not=l'\leq n+2}(1+\eps_l^2)(\hR_{l,l'}-q)^2.
$$
On the other hand, when $l=l'$ we get directly
\begin{equation}
\mbox{Cov}\Bigl(
H_t(\vsi^l),\frac{\partial}{\partial t} H_t(\vsi^{l})
\Bigr)
=\frac{\beta^2}{2}\Bigl(\xi(1) -\eps_l^2 \xi'(q)\Bigr).
\label{cov2}
\end{equation}
Next, we simplify the third term on the right hand side of (\ref{Hamintder}).
(\ref{xdel}) implies
$$
\Bigl|
\sqrt{\frac{N-\eps^2}{N-1}}
-\Bigl(1+\frac{1-\eps^2}{2(N-1)}\Bigr)
\Bigr|\leq L \frac{(1-\eps^2)^2}{(N-1)^2}
$$
and, therefore,
$$
(N-1)\Bigl(\sqrt{\frac{N-\eps^2}{N-1}}-1\Bigr)
-\frac{1-\eps^2}{2} =\R_1.
$$
We can write
\begin{eqnarray}
h\sum_{i\leq N-1}\hs_i^l\Bigl(\sqrt{\frac{N-\eps_l^2}{N-1}}-1\Bigr)
&=& 
\frac{h}{2}\hR_l (1-\eps_l^2) + h\R_1
\nonumber
\\
&=&
\frac{h}{2}a_l(\hR_l -r) 
+ \frac{hr}{2}(1-\eps_l^2)+h\R_1,
\label{hsimple}
\end{eqnarray}
where in the last line we used the definition of $a_l$ in (\ref{ais}).
Finally, using (\ref{cov1}), (\ref{cov2}) and (\ref{hsimple}), 
Gaussian integration by parts in (\ref{der1}) gives,
$$
\nu_t'(f)= \mbox{I}+\mbox{II}+\mbox{III}+\mbox{IV}+\mbox{V}+\mbox{VI}+\nu_t(f\R),
$$
where I is created by the first term in (\ref{hsimple}):
$$
\mbox{\rm I}=\frac{h}{2}\sum_{l\leq n} \nu_t(f a_l(\hR_l -r)) -
n\frac{h}{2}\nu_t(f a_{n+1}(\hR_{n+1}-r)),
$$
II is created by the first term in (\ref{cov1}):
\begin{eqnarray*}
\mbox{\rm II} &=&
\beta^2 \sum_{1\leq l\not=l'\leq n} \nu_t(f a_{l,l'}(\hR_{l,l'}-q))
-2n\beta^2 \sum_{l\leq n} \nu_t(f a_{l,n+1}(\hR_{l,n+1}-q))
\nonumber
\\
&+&
n(n+1)\beta^2 \nu_t(f a_{n+1,n+2}(\hR_{n+1,n+2}-q)),
\end{eqnarray*}
III is created by the second term in (\ref{hsimple}):
$$
\mbox{\rm III}=-\frac{hr}{2}\Bigl(\sum_{l\leq n} \nu_t(f\eps_l^2) 
- n\nu_t(f \eps_{n+1}^2)\Bigr),
$$
IV is created by the second term in (\ref{cov1}):
\begin{eqnarray*}
\mbox{\rm IV} &=&
-\frac{\beta^2}{4}q\xi'(q)\Bigl( 
\sum_{1\leq l\not=l'\leq n} \nu_t(f (\eps_l^2+\eps_{l'}^2))
-2n\sum_{l\leq n} \nu_t(f (\eps_l^2+\eps_{n+1}^2))
\nonumber
\\
&+&
n(n+1)\nu_t(f (\eps_{n+1}^2+\eps_{n+2}^2))\Bigr),
\end{eqnarray*}
V is created by (\ref{cov2}):
$$
\mbox{\rm V}=-\frac{\beta^2}{2}\xi'(q)\Bigl(\sum_{l\leq n} \nu_t(f\eps_l^2) 
- n\nu_t(f \eps_{n+1}^2)\Bigr),
$$
and VI is created by the last term in (\ref{Hamintder}):
$$
\mbox{\rm VI}=\frac{1}{2}b\Bigl(\sum_{l\leq n} \nu_t(f\eps_l^2) 
- n\nu_t(f \eps_{n+1}^2)\Bigr).
$$
Using that by symmetry, $\nu_t(f \eps_{n+1}^2) = \nu_t(f\eps_{n+2}^2),$
and counting terms in IV it is easy to see that
$$
\mbox{\rm IV}=\frac{\beta^2}{2}q\xi'(q)\Bigl(\sum_{l\leq n} \nu_t(f\eps_l^2) 
- n\nu_t(f \eps_{n+1}^2)\Bigr).
$$
Since, by definition, $b=hr +\beta^2(1-q)\xi'(q),$ we have
III+IV+V+VI$=0.$
This finishes the proof of Theorem \ref{Thder}.
\qed

The goal of the above interpolation is to relate $\nu(f)$ to $\nu_0(f)$
because for proper choices of the function $f$ one can compute 
(or accurately estimate) $\nu_0(f)$ due to the special form of 
the Hamiltonian (\ref{Hamzero}) at $t=0.$ Therefore, in order for
this interpolation to be useful, the derivative (\ref{der}) should
be small. This fact is contained in the following two results.

\begin{theorem}\label{epscontrol}
If $\beta$ and $h$ are small enough, we can find a constant $L>0$
such that
\begin{equation}
\nu_t\Bigl(\exp \frac{1}{L}\eps^2\Bigr) \leq L
\label{expeps}
\end{equation}
for all $t\in[0,1].$
\end{theorem}

\begin{theorem}\label{Big}
If $\beta$ and $h$ are small enough then for any $K>0$ we can find
$L>0$ such that
\begin{eqnarray}
&&
\nu_t\Bigl(I\Bigl(|\hR_{1,2}-q|\geq L\Bigl(\frac{\log N}{N}\Bigr)^{1/4}\Bigr)\Bigr)
\leq \frac{L}{N^{K}},
\label{Big1}
\\
&&
\nu_t\Bigl(I\Bigl(|\hR_{1}-r|\geq L\Bigl(\frac{\log N}{N}\Bigr)^{1/4}\Bigr)\Bigr)
\leq \frac{L}{N^{K}}
\label{Big2}
\end{eqnarray}
for all $t\in [0,1].$
\end{theorem}
We will prove Theorem \ref{epscontrol} in Section \ref{Seceps}
and Theorem \ref{Big} in Section \ref{SecR}.
It is rather clear that they will provide the necessary control of 
each term in the derivative (\ref{der}), which will be demonstrated
in the next section.

Next we will explain what happens at the end of the interpolation at $t=0.$
Let us start by writing the integration over $S_N$ as a double integral
over $\eps$ and $(\sigma_1,\ldots,\sigma_{N-1}).$ Let $\lambda_{N}^{\rho}$ denote
the area measure on the sphere $S_N^{\rho}$ of radius $\rho$ in $\Reals^N,$
and let $|S_N^{\rho}|$ denote its area, i.e. 
$|S_N^{\rho}|=\lambda_{N}^{\rho}(S_N^{\rho}).$
 Then,
\begin{eqnarray}
&&
\int\limits_{S_N}f(\vsi)d\lambda_N(\vsi)=\frac{1}{|S_N^{\sqrt{N}}|}
\int\limits_{S_N^{\sqrt{N}}}f(\sigma_1,\ldots,\sigma_N)d\lambda_N^{\sqrt{N}}(\sigma_1,\ldots,\sigma_N)
\nonumber
\\
&&
=
\frac{1}{|S_N^{\sqrt{N}}|}
\int\limits_{-\sqrt{N}}^{\sqrt{N}}\frac{d\eps}{\sqrt{1-\eps^2/N}}
\int\limits_{S_{N-1}^{\sqrt{N-\eps^2}}} f(\sigma_1,\ldots,\sigma_{N-1},\eps)
d\lambda_{N-1}^{\sqrt{N-\eps^2}}(\sigma_1,\ldots,\sigma_{N-1})
\nonumber
\\
&&
=
\int\limits_{-\sqrt{N}}^{\sqrt{N}}\frac{|S_{N-1}^{\sqrt{N-\eps^2}}|}{|S_N^{\sqrt{N}}|}
\frac{d\eps}{\sqrt{1-\eps^2/N}}
\int\limits_{S_{N-1}} f\Bigl(\hs_1\sqrt{\frac{N-\eps^2}{N-1}},\ldots,\hs_{N-1}
\sqrt{\frac{N-\eps^2}{N-1}},\eps\Bigr)
d\lambda_{N-1}(\hS)
\nonumber
\\
&&
=a_N
\int\limits_{-\sqrt{N}}^{\sqrt{N}}
d\eps\Bigl(1-\frac{\eps^2}{N}\Bigr)^{\frac{N-3}{2}}
\int\limits_{S_{N-1}} f\Bigl(\hS\sqrt{\frac{N-\eps^2}{N-1}},\eps\Bigr)
d\lambda_{N-1}(\hS),
\label{Fub1}
\end{eqnarray}
where $a_N=|S_{N-1}^1|/(|S_{N}^1|\sqrt{N})\to (2\pi)^{-1/2}$ as can be seen
by taking $f=1.$
In particular, if 
$$
f(\vsi)=f_1(\eps)f_2(\hS)
$$ 
then
\begin{equation} 
\int\limits_{S_N}
f(\vsi) d\lambda_N(\vsi)
=a_N 
\int\limits_{-\sqrt{N}}^{\sqrt{N}}
f_1(\eps)\Bigl(1-\frac{\eps^2}{N}\Bigr)^{\frac{N-3}{2}}d\eps
\int\limits_{S_{N-1}} f_2(\hS) d\lambda_{N-1}(\hS).
\label{Fub2}
\end{equation}
Since the Hamiltonian (\ref{Hamzero}) decomposed into the sum of terms 
that depend only on $\eps$ or only on $\hS,$ (\ref{Fub2}) implies
that
\begin{equation}
\la f\ra_0 = \la f_1\ra_0 \la f_2\ra_0
\label{decop}
\end{equation}
where
\begin{equation}
\bigl\la f_1(\eps)\bigr\ra_0 = \frac{1}{Z_1}\int\limits_{-\sqrt{N}}^{\sqrt{N}}
f_1(\eps)\Bigl(1-\frac{\eps^2}{N}\Bigr)^{\frac{N-3}{2}}
\exp\Bigl(a\eps-\frac{1}{2}b\eps^2\Bigr) d\eps,
\label{decopone}
\end{equation}
$$
Z_1= \int\limits_{-\sqrt{N}}^{\sqrt{N}}
\Bigl(1-\frac{\eps^2}{N}\Bigr)^{\frac{N-3}{2}}
\exp\Bigl(a\eps-\frac{1}{2}b\eps^2\Bigr) d\eps,
$$
and
\begin{equation}
\bigl\la f_2(\hS)\bigr\ra_0 = \frac{1}{Z_2}\int\limits_{S_{N-1}}
f_2(\hS)\exp\Bigl(H_{N-1}(\hS)+h\sum_{i\leq N-1}\hs_i\Bigr) d\lambda_{N-1}(\hS),
\label{decoptwo}
\end{equation}
$$
Z_2= 
\int\limits_{S_{N-1}}
\exp\Bigl(H_{N-1}(\hS)+h\sum_{i\leq N-1}\hs_i\Bigr) d\lambda_{N-1}(\hS).
$$
Using (\ref{decop}), (\ref{decopone}), we will be able to compute the moments
$\nu_0(\eps_1^{k_1}\ldots \eps_n^{k_n})$ for integer $k_i\geq 0,$
which is an important part of the second moment computations and
of the cavity method in general. This is done as follows.
Let us recall (\ref{b}), (\ref{a}) and define $\gamma_0=1, 
\gamma_1=a/(b+1)$ and, recursively, for $k\geq 2$
\begin{equation}
\gamma_k=\frac{a}{b+1}\gamma_{k-1} +\frac{k-1}{b+1}\gamma_{k-2}.
\label{gammas}
\end{equation}
The following Theorem holds.
\begin{theorem}\label{zeromoments}
For small enough $\beta>0,$
\begin{equation}
\Bigl|
\nu_0(\eps_1^{k_1}\ldots \eps_n^{k_n})
-\e \gamma_{k_1}\ldots \gamma_{k_n}
\Bigr|\leq \frac{L}{N},
\end{equation}
where a constant $L$ is independent of $N.$
\end{theorem}
This Theorem will be proved in Section \ref{Seceps} below.

\section{Second moment computations.}

Let us introduce the following seven functions
\begin{eqnarray}
&&
f_1=(R_{1,2}-q)^2,\,\,f_2=(R_{1,2}-q)(R_{1,3}-q),\,\,
f_3=(R_{1,2}-q)(R_{3,4}-q)
\label{fs}
\\
&&
f_4=(R_{1,2}-q)(R_1-r),\,\,f_5=(R_{1,2}-q)(R_3-r),\,\,
f_6=(R_1-r)^2,\,\, f_7=(R_1-r)(R_2-r)
\nonumber
\end{eqnarray}
and let $\vec{v}_N=(\nu(f_1),\ldots,\nu(f_7)).$ In this section we
will compute a vector $N\vec{v}_N$ up to the terms of order $o(1).$
As we mentioned above, it is likely that with more effort one can
prove the central limit theorem for the joint distribution of
$$
\sqrt{N}(R_{1,2}-q),\,\,\, 
\sqrt{N}(R_{1,3}-q),\,\,\, 
\sqrt{N}(R_{3,4}-q),\,\,\,
\sqrt{N}(R_{1}-r),\,\,\,
\sqrt{N}(R_{2}-r),\,\,\,   
$$
so the computation of this section identifies the covariance
matrix of the limiting Gaussian distribution.
To describe our main result let us first summarize several computations 
based on Theorem \ref{zeromoments}. 
The definition (\ref{gammas}) implies that
\begin{equation}
\gamma_1=\frac{a}{b+1},\,\,
\gamma_2=\Bigl(\frac{a}{b+1}\Bigr)^2+\frac{1}{b+1},\,\,
\gamma_3=\Bigl(\frac{a}{b+1}\Bigr)^3+\frac{3a}{(b+1)^2}.
\label{gams}
\end{equation}
The definition (\ref{b}) and (\ref{critq}) imply that
$$
\frac{1}{b+1}=\frac{1}{1+(1-q)(\beta^2 \xi'(q) +h^2)}=1-q.
$$
Therefore,
\begin{equation}
\e\frac{a}{b+1}=(1-q)\e a=(1-q)h=r,\,\,
\end{equation}
\begin{equation}
\e \Bigl(\frac{a}{b+1}\Bigr)^2 =(1-q)^2\e a^2
=(1-q)^2 (\beta^2\xi'(q)+h^2)=q
\end{equation}
where we used (\ref{critq}) again, and
\begin{equation}
W:=\e \Bigl(\frac{a}{b+1}\Bigr)^3 =(1-q)^3\e a^3
=(1-q)^3 (3\beta^2 \xi'(q) h+h^3),
\end{equation}
\begin{equation}
U:=\e \Bigl(\frac{a}{b+1}\Bigr)^4 =(1-q)^4\e a^4
=(1-q)^4(h^4+ 6\beta^2 h^2 \xi'(q)+3\beta^4 \xi'(q)^2). 
\label{through}
\end{equation}
For simplicity of notations let us write
$$
x\sim y \,\,\,\mbox{ if }\,\,\,
x=y+O\bigl(N^{-1}\bigr).
$$
Then it is trivial to check that Theorem \ref{zeromoments} and 
(\ref{gams}) - (\ref{through}) imply the following relations:
\begin{eqnarray}
&&
\nu_0(\eps_1)\sim r,\,\,\,
\nu_0(\eps_1\eps_2)\sim q,\,\,\,
\nu_0(\eps_1^2)\sim 1,\,\,\,
\nu_0(\eps_1\eps_2\eps_3)\sim  W,\,\,\,
\nonumber
\\
&&
\nu_0(\eps_1\eps_2^2)\sim W+h(1-q)^2,\,\,\,
\nu_0(\eps_1^3)\sim W+3h(1-q)^2,
\nonumber
\\
&&
\nu_0(\eps_1^2\eps_2^2)\sim U+1-q^2,\,\,\,
\nu_0(\eps_1\eps_2\eps_3^2)\sim U+q-q^2
\nonumber
\\
&&
\nu_0(\eps_1\eps_2^3)\sim U+3q-3q^2,\,\,\,\,
\nu_0(\eps_1\eps_2\eps_3\eps_4)\sim U.
\label{relations}
\end{eqnarray}
Let us recall the definitions
$a_l$ and $a_{l,l'}$ in (\ref{ais}). Using relations (\ref{relations})
it is now straightforward to compute the following 
nine quantities
\begin{eqnarray}
&&
\nu_0(a_{1,2}(\eps_1\eps_2-q))\sim Y_1,\,\,
\nu_0(a_{1,3}(\eps_1\eps_2-q))\sim Y_2,\,\,
\nu_0(a_{3,4}(\eps_1\eps_2-q))\sim Y_3,
\nonumber
\\
&&
\nu_0(a_1(\eps_1\eps_2-q))\sim Y_4,\,\,
\nu_0(a_3(\eps_1\eps_2-q))\sim Y_5,\,\,
\nu_0(a_{1,2}(\eps_1-r))\sim Y_6,
\nonumber
\\
&&
\nu_0(a_{2,3}(\eps_1-r))\sim Y_7,\,\,
\nu_0(a_1(\eps_1-r))\sim Y_8,\,\,
\nu_0(a_2(\eps_1-r))\sim Y_9,
\label{relations2}
\end{eqnarray}
where $Y_1,\ldots,Y_9$ are functions of $q,r,h, U,W.$
We omit the explicit formulas for $Y_j$s since they do not serve
any particular purpose in the sequel.
Let us define a 7$\times$7 matrix $M$ that consists of four
blocks
\begin{equation}
M=\left(
\begin{array}{cc}
M_1 & O_1 \\
O_2 & M_2
\end{array}
\right)
\label{M}
\end{equation}
where $O_1$ is a 3$\times$2 matrix and $O_2$ is a 4$\times$3 matrix 
both entirely consisting of zeros, 
$$
M_1=\left(
\begin{array}{ccccc}
2\beta^2 Y_1 & -8\beta^2 Y_2 & 6\beta^2 Y_3 & h Y_4 & -h Y_5
\\
2\beta^2 Y_2 & 2\beta^2(Y_1-2 Y_2-3Y_3) & 6\beta^2(-Y_2+2Y_3) &
\frac{h}{2}(Y_4+Y_5) & \frac{h}{2}(Y_4-3Y_5)
\\
2\beta^2 Y_3 & 8\beta^2(Y_2-2Y_3) & 2\beta^2(Y_1-8Y_2+10 Y_3) &
h Y_5 & h (Y_4- 2Y_5)
\end{array}
\right),
$$
$$
M_2=\left(
\begin{array}{cccc}
2\beta^2(Y_1-2 Y_2) & 2\beta^2(-2Y_2+ 3Y_3) &
(h/2)Y_4 & (h/2)(Y_4-2Y_5)
\\
2\beta^2(2Y_2-3 Y_4) & 2\beta^2(Y_1-6Y_2+6 Y_3) &
(h/2)Y_5 & (h/2)(2Y_4-3Y_5)
\\
-2\beta^2Y_6 & 2\beta^2 Y_7 & (h/2)Y_8 & -(h/2)Y_9
\\
2\beta^2(Y_6-2Y_7) & 2\beta^2(-2Y_6+3Y_7) &
(h/2)Y_9 & (h/2)(Y_8-2Y_9)
\end{array}
\right).
$$
Finally, we define a vector $\vec{v}=(v_1,\ldots,v_7)$ by
\begin{eqnarray}
&&
v_1=(1-q)U+1-4q^2+3q^3,\,\,\, 
v_2=(1-q)U+q(1-q)(1-2q)
\nonumber
\\
&&
v_3=(1-q)U-q^2(1-q),\,\,\,
v_4=W-\frac{1}{2}rU+\frac{1}{2}r(2-6q+3q^2)
\nonumber
\\
&&
v_5=W-\frac{1}{2}rU+\frac{1}{2}r(-2q+q^2),\,\,\,
v_6=-\frac{1}{2}rW+1+\frac{1}{2}r^2(-4+3q)
\nonumber
\\
&&
v_7=-\frac{1}{2}rW+q+\frac{1}{2}r^2(-2+q).
\label{vees}
\end{eqnarray}
We are now ready to formulate the main result of this section.
\begin{theorem}\label{Main}
For small enough $\beta$ and $h$ we have
\begin{equation}
(I-M)\vec{v}_N^T = \frac{1}{N} \vec{v}^T + o\bigl(N^{-1}\bigr).
\label{secmom}
\end{equation}
\end{theorem}
Here $\vec{v}^T$ denotes the transpose of vector $\vec{v}.$
Notice that each entry in the matrix $M$ has either a factor of $\beta^2$
or $h$ and, therefore, for small enough $\beta$ and $h$ the matrix $(I-M)$
will be invertible, in which case Theorem \ref{Main} implies
$$
\vec{v}_N^T=\frac{1}{N}(I-M)^{-1}\vec{v}^T +o(N^{-1}).
$$
In the remainder of this section we will prove Theorem \ref{Main}.

For each function $f_l$ in (\ref{fs}), we will define $\hat{f}_l$
by replacing each occurrence of $R$ by $\hR,$ i.e. $\hat{f_1}=(\hR_{1,2}-q)^2,$
$\hat{f}_2 = (\hR_{1,2}-q)(\hR_{1,3}-q)$ etc. Next, we introduce functions
\begin{eqnarray}
&&
f_1'=(\eps_1\eps_2-q)(R_{1,2}-q),\,\,f_2'=(\eps_1\eps_2-q)(R_{1,3}-q),\,\,
f_3'=(\eps_1\eps_2-q)(R_{3,4}-q)
\nonumber
\\
&&
f_4'=(\eps_1\eps_2-q)(R_1-r),\,\,f_5'=(\eps_1\eps_2-q)(R_3-r),
f_6'=(\eps_1-r)(R_1-r),
\nonumber
\\
&&
f_7'=(\eps_1-r)(R_2-r).
\label{fprimes}
\end{eqnarray}
As in the classical cavity method in \cite{SG},
we introduce these functions because, first of all, by symmetry,
\begin{equation}
\nu(f_l)=\nu(f_l')
\label{symme} 
\end{equation}
and, second of all, emphasizing the last coordinate in $f_l'$ 
is perfectly suited for the application of the cavity method.
As above, for each function $f_l'$ we will define $\hat{f}_l'$
by replacing each occurrence of $R$ by $\hR,$ i.e. 
$\hat{f}_1'=(\eps_1\eps_2-q)(\hR_{1,2}-q)$ etc.

To simplify the notations we will write $x\approx y$ whenever
\begin{equation}
|x-y|=o\Bigl(
\frac{1}{N}+\nu_0((\hR_{1,2}-q)^2)+\nu_0((\hR_1-r)^2)
\Bigr).
\label{equi}
\end{equation}
The proof of Theorem \ref{Main} will be based on the following.
\begin{theorem}\label{Theq}
For small enough $\beta$ and $h,$ for all $l\leq 7,$
\begin{equation}
\nu_0(\hat{f}_l)\approx \nu_0(f_l')+\nu_0'(\hat{f}_l').
\label{Theqeq}
\end{equation}
\end{theorem}
We will start with a couple of lemmas.
\begin{lemma}
If $f\geq 0$ and $\|f\|_{\infty}$ is bounded independently of $N$ 
then for any $K>0$ we can find $L>0$ such that
\begin{equation}
\nu_t(f)\leq L\Bigl(N^{-K} +\nu_0(f)\Bigr).
\label{derzero}
\end{equation}
\end{lemma}
\textbf{Proof.}
The derivative $\nu_t'(f)$ in (\ref{der}) consists of a finite sum of terms
of the type $\nu_t(f\pepsi g)$ where $\pepsi$ is some polynomial in the
last coordinates $(\eps_l)$ and $g$ is one of the following:
\begin{equation}
\hR_{l,l'}-q,\,\,\, \hR_l-r,\,\,\, (\hR_{l,l'}-q)^2, \,\,\, N^{-1}. 
\label{hs}
\end{equation}
Theorem \ref{epscontrol} and Chebyshev's inequality imply
$$
\nu_t\Bigl(I\bigl(|\eps_l|\geq \log N\bigr)\Bigr)\leq LN^{-K}
$$
and combining this with Theorem \ref{Big} yields that for any 
$g$ in (\ref{hs}),
$$
\nu_t\Bigl(I\bigl(|\pepsi g|\geq N^{-1/8}\bigr)\Bigr)\leq LN^{-K}.
$$
Therefore, one can control the derivative
\begin{equation}
|\nu_t'(f)|\leq LN^{-K}+LN^{-1/8}\nu_t(f) \leq L(N^{-K}+\nu_t(f))
\label{dereps}
\end{equation}
and (\ref{derzero}) follows by integration.
\qed

\begin{lemma}\label{hatzero}
For small enough $\beta$ and $h$ and all $l\leq 7$ we have
\begin{equation}
\nu(f_l)\approx\nu_0(\hat{f}_l) \,\,\mbox{ and }\,\,
\nu_0'(f_l')\approx \nu_0'(\hat{f}_l').
\label{hatzeroeq}
\end{equation}
\end{lemma}
\textbf{Proof.}
We will only consider the case $l=1,$ $f_1=(R_{1,2}-q)^2$,
since other cases are similar. We have
\begin{eqnarray}
&&
\Bigl|
\nu\bigl((R_{1,2}-q)^2\bigr)-\nu_0\bigl((R_{1,2}-q)^2\bigr)
\Bigr|
\leq \sup_{t}\Bigl(\nu_t'\bigl((R_{1,2}-q)^2\bigr)\Bigr)
\label{ablem}
\\
&&
\leq \sup_{t}\Bigl(
LN^{-K}+LN^{-1/8}\nu_t\bigl((R_{1,2}-q)^2\bigr)\Bigr)
\leq
\Bigl(
LN^{-K}+LN^{-1/8}\nu_0\bigl((R_{1,2}-q)^2\bigr)
\Bigr)
\nonumber
\end{eqnarray}
where in the second line we used (\ref{dereps}) and then (\ref{derzero}).
Since by (\ref{ars})
\begin{equation}
R_{1,2}-q=(\hR_{1,2}-q)
+\Bigl(\sqrt{\Bigl(1-\frac{\eps_1^2}{N}\Bigr)\Bigl(1-\frac{\eps_2^2}{N}\Bigr)} -1\Bigr)
\hR_{1,2}
+\frac{1}{N}\eps_1\eps_2,
\label{arsq}
\end{equation}
squaring both sides and using (\ref{epsils}) yields
$$
\Bigl|
(R_{1,2}-q)^2 -(\hR_{1,2}-q)^2
\Bigr| \leq
\frac{1}{N}\pepsi|\hR_{1,2}-q|+\frac{1}{N^2}\pepsi,
$$
where from now on $p_{\eps}$ denotes a quantity such that
$$
|p_{\eps}|\leq L(1+\sum_{l} \eps_l^4).
$$
Therefore,
$$
\Bigl|
\nu_0\bigl((R_{1,2}-q)^2\bigr)-\nu_0\bigl((\hR_{1,2}-q)^2\bigr)
\Bigr|\leq \frac{1}{N}\nu_0\bigl(\pepsi|\hR_{1,2}-q|\bigr)+\frac{1}{N^2}\nu_0(\pepsi)
=o(N^{-1})
$$
by Theorems \ref{epscontrol} and \ref{Big}.
Thus, (\ref{ablem}), implies the first part of (\ref{hatzeroeq}).
To prove the second part of (\ref{hatzeroeq}) we notice that
$$
\Bigl|
f_1'-\hat{f}_1'
\Bigr|=
\Bigl|
(\eps_1\eps_2-q)(R_{1,2}-\hR_{1,2})
\Bigr| \leq \frac{1}{N}\pepsi
$$
by (\ref{ars}) and (\ref{epsils}). Since each term in the derivatives
$\nu_0'(f_1')$ and $\nu_0'(\hat{f}_1')$ will contain another factor from 
the list (\ref{hs}), Theorems \ref{epscontrol} and \ref{Big} imply the result.
\qed

\textbf{Proof of Theorem \ref{Theq}.}
We start by writing
$$
\Bigl|
\nu(f_l')-\nu_0(f_l')-\nu_0'(f_l')
\Bigr|
\leq \sup_{t}\Bigl|\nu_t''(f_l')\Bigr|.
$$
If we can show that 
\begin{equation}
\sup_{t}\Bigl|\nu_t''(f_l')\Bigr| \approx 0
\label{secder}
\end{equation}
and, thus,
$
\nu(f_l')\approx \nu_0(f_l')+\nu_0'(f_l'),
$
then Lemma \ref{hatzero} and (\ref{symme}) will imply
$$
\nu_0(\hat{f}_l)\approx \nu(f_l)=\nu(f_l')
\approx \nu_0(f_l')+\nu_0'(f_l')
\approx \nu_0(f_l')+\nu_0'(\hat{f}_l'),
$$
which is precisely the statement of Theorem \ref{Theq}.
To prove (\ref{secder}) we note that by (\ref{der})
the second derivative $\nu_t''(f_l')$ will consist of the
finite sum of terms of the type $f_l'\pepsi g_1 g_2$
where $g_1, g_2$ are from the list (\ref{hs}). Clearly,
$$
|g_1g_2|\leq L\Bigl(
\frac{1}{N^2}+(\hR_{l,l'}-q)^2+(\hR_{l''}-r)^2
\Bigr)
$$
and since each $f_l'$ contain another small factor 
$(R_{l,l'}-q)$ or $(R_l-r),$ Theorems \ref{epscontrol}
and \ref{Big} imply (\ref{secder}).
\qed

We are now ready to prove Theorem \ref{Main}.

\textbf{Proof of Theorem \ref{Main}.}
Let us first note that $\nu_0(\hat{f}_l)$ is defined
exactly the same way as $\nu(f_l)$ for $N-1$ instead of $N.$
In other words, 
$$
\vec{v}_N^0 :=(\nu_0(\hat{f}_1),\ldots,\nu_0(\hat{f}_7))=\vec{v}_{N-1}
$$ 
and, therefore, it is enough to prove that
\begin{equation}
(I-M){\vec{v}_N^0}^T=\frac{1}{N}\vec{v}^T+o(N^{-1}).
\label{IM}
\end{equation}
Replacing $1/N$ by $1/(N-1)$ on the right hand side is not necessary since
the difference is of order $N^{-2}.$ Each equation in the system of 
equations (\ref{IM}) will follow from the corresponding equation (\ref{Theqeq}).
Namely, we will show that 
\begin{equation}
(\nu_0(f_1'),\ldots,\nu_0(f_7'))\approx \frac{1}{N}\vec{v}\,\,\,
\mbox{ and } \,\,\,
(\nu_0'(\hat{f}_1'),\ldots,\nu_0'(\hat{f}_7'))^T\approx M{\vec{v}_{N}^0}^T.
\label{show}
\end{equation}
Then (\ref{Theqeq}) will imply that ${\vec{v}_{N}^0}^T\approx N^{-1}\vec{v}^T
+M{\vec{v}_N^0}^T.$ However, since the definition (\ref{equi}) means that the error
in each equation is of order $o(N^{-1}+\nu_0(\hat{f}_1)+\nu_0(\hat{f}_6)),$
this system of equation can be rewritten as
$$
(I-M-\E_N){\vec{v}_N^0}^T=\frac{1}{N}\vec{v},
$$
where the matrix $\E_N$ is such that $\|\E_N\|=o(1).$
Therefore, whenever the matrix $I-M$ is invertible 
(for example, for small $\beta$ and $h$) we have for $N$ large enough
$$
{\vec{v}_N^0}^T=\frac{1}{N}(I-M-\E_N)^{-1}\vec{v}^T
=\frac{1}{N}(I-M)^{-1}\vec{v}^T +o(N^{-1}).
$$
Hence, to finish the proof we need to show (\ref{show}).
We will only carry out the computations for $l=1$ since all
other cases are similar. Let us start by proving that
$\nu_0\bigl((\eps_1\eps_2-q)(R_{1,2}-q)\bigr)\approx v_1.$ 
Using (\ref{arsq}) and (\ref{decop}), we write 
\begin{eqnarray*}
\nu_0\bigl((\eps_1\eps_2-q)(R_{1,2}-q)\bigr)
&=&
\frac{1}{N}\nu_0\bigl(\eps_1\eps_2(\eps_1\eps_2-q)\bigr)
+\nu_0\bigl((\eps_1\eps_2-q)\bigr)\nu_0\bigl((\hR_{1,2}-q)\bigr)
\\
&+&
q\nu_0\Bigl(
(\eps_1\eps_2 -q)\Bigl(
\sqrt{\Bigl(1-\frac{\eps_1^2}{N}\Bigr)\Bigl(1-\frac{\eps_2^2}{N}\Bigr)}-1
\Bigr)
\Bigr)
\\
&+&
\nu_0\Bigl(
(\eps_1\eps_2 -q)\Bigl(
\sqrt{\Bigl(1-\frac{\eps_1^2}{N}\Bigr)\Bigl(1-\frac{\eps_2^2}{N}\Bigr)}-1
\Bigr)
\Bigr)
\nu_0\bigl(\hR_{1,2}-q\bigr).
\end{eqnarray*}
Using (\ref{epsils}), one can bound the last term by
$$
\frac{1}{N}\nu_0(\pepsi)\Bigl|\nu_0(\hR_{1,2}-q)\Bigr|=o(N^{-1}),
$$
by Theorems \ref{epscontrol} and \ref{Big}.
The term
$$
\nu_0\bigl((\eps_1\eps_2-q)\bigr)\nu_0\bigl((\hR_{1,2}-q)\bigr)
=o(N^{-1})
$$
by Theorem \ref{Big} and the second relation in (\ref{relations}),
i.e. $\nu_0(\eps_1\eps_2-q)\sim 0.$ Finally, we use
$$
\Bigl|
\sqrt{\Bigl(1-\frac{\eps_1^2}{N}\Bigr)\Bigl(1-\frac{\eps_2^2}{N}\Bigr)}-1
+\frac{\eps_1^2}{2N}+\frac{\eps_2^2}{2N}
\Bigr|\leq \frac{1}{N^2}\pepsi
$$
to observe that
\begin{eqnarray*}
q\nu_0\Bigl(
(\eps_1\eps_2 -q)\Bigl(
\sqrt{\Bigl(1-\frac{\eps_1^2}{N}\Bigr)\Bigl(1-\frac{\eps_2^2}{N}\Bigr)}-1
\Bigr)
\Bigr)
&\approx&
-\frac{1}{2N}q\nu_0\bigl((\eps_1\eps_2-q)(\eps_1^2+\eps_2^2)\bigr)
\\
&=&
-\frac{1}{N}q\nu_0\bigl((\eps_1\eps_2-q)\eps_1^2\bigr)
\end{eqnarray*}
by symmetry and, therefore,
\begin{eqnarray*}
\nu_0\bigl((\eps_1\eps_2-q)(R_{1,2}-q)\bigr)
&\approx&
\frac{1}{N}\Bigl(
\nu_0\bigl(\eps_1\eps_2(\eps_1\eps_2-q)\bigr)
-q\nu_0\bigl((\eps_1\eps_2-q)\eps_1^2\bigr)
\Bigr)
\\
&=&
\frac{1}{N}\Bigl(
\nu_0\bigl(\eps_1^2\eps_2^2\bigr)-q\nu_0\bigl(\eps_1\eps_2\bigr)
-q\nu_0\bigl(\eps_1^3\eps_2\bigr)+q^2\nu_0\bigl(\eps_1^2\bigr)
\Bigr)
\approx v_1
\end{eqnarray*}
by using (\ref{relations}) and comparing with the definition
of $v_1$ in (\ref{vees}).

Next, we need to show the second part of (\ref{show}) for
$l=1,$ i.e. 
$$
\nu_0'((\eps_1\eps_2-q)(\hR_{1,2}-q))\approx 
\bigl(M{\vec{v}_N^0}^T\Bigr)_1.
$$ 
We use (\ref{der}) for $n=2$ to write $\nu_0'((\eps_1\eps_2-q)(\hR_{1,2}-q))$ as
\begin{eqnarray*}
&&
h\nu_0(a_1(\eps_1\eps_2-q)(\hR_{1,2}-q)(\hR_{1}-r))
-h\nu_0(a_3(\eps_1\eps_2-q)(\hR_{1,2}-q)(\hR_{3}-r))
\\
&+&
2\beta^2h\nu_0(a_{1,2}(\eps_1\eps_2-q)(\hR_{1,2}-q)^2)
-8\beta^2 h\nu_0(a_{1,3}(\eps_1\eps_2-q)(\hR_{1,2}-q)(\hR_{1,3}-q))
\\
&+&
6\beta^2h\nu_0(a_{3,4}(\eps_1\eps_2-q)(\hR_{1,2}-q)(\hR_{3,4}-q))
+\nu_0((\eps_1\eps_2-q)(\hR_{1,2}-q)\R)
\\
&\approx &
2\beta^2 Y_1 \nu_0(\hat{f}_1) -8\beta^2 Y_2\nu_0(\hat{f}_2)
+6\beta^2 Y_3\nu_0(\hat{f}_3)+h Y_4 \nu_0(\hat{f}_4)
-h Y_5\nu_0(\hat{f}_5)+\nu_0(\hat{f}_1'\R)
\\
&=&
\bigl(M{\vec{v}_N^0}^T\Bigr)_1 + \nu_0(\hat{f}_1'\R),
\end{eqnarray*}
where in second to last line we used (\ref{relations2}) and the
last line follows by comparison with the definition of $M$
in (\ref{M}). Finally, since clearly $\nu_0(\hat{f}_1'\R)\approx 0$
by Theorems \ref{epscontrol} and \ref{Big}, this finishes the
proof of Theorem \ref{Main}
\qed

\section{Control of the last coordinate.}\label{Seceps}

In this section we will prove Theorems \ref{epscontrol} and \ref{zeromoments}.
We start with the following.
\begin{lemma}\label{epszero}
If $c_0<1$ then for $\beta$ small enough, 
$$
\nu_0\bigl(\exp c_0\eps^2\bigr)\leq L.
$$
\end{lemma}
\textbf{Proof.}
By (\ref{decop}) and using $1-x\leq \exp(-x)$,
\begin{eqnarray}
\bigl\la\exp c_0\eps^2\bigr\ra_0
&=&
\frac{1}{Z_1}\int\limits_{-\sqrt{N}}^{\sqrt{N}}
\Bigl(1-\frac{\eps^2}{N}\Bigr)^{\frac{N-3}{2}}
\exp\Bigl(a\eps-\frac{1}{2}(b-c_0)\eps^2\Bigr) d\eps,
\label{ceps}
\\
&\leq&
\frac{1}{Z_1}\int\limits_{-\sqrt{N}}^{\sqrt{N}}
\exp\Bigl(a\eps-\frac{1}{2}(b-c_0+1-3N^{-1})\eps^2\Bigr) d\eps
\leq
\frac{1}{Z_1}L\exp(La^2)
\nonumber
\end{eqnarray}
since for $c_0<1$ we have $b+1-3N^{-1}-c_0>0$
for large enough $N.$
On the other hand, one can show that
\begin{equation}
Z_1\geq \frac{1}{L}\exp(-La^2).
\label{bound1}
\end{equation}
Indeed, using that $1-x\geq \exp(-Lx)$ for $x\leq 1/2,$
\begin{eqnarray*}
Z_1
&=& 
\int\limits_{-\sqrt{N}}^{\sqrt{N}}
\Bigl(1-\frac{\eps^2}{N}\Bigr)^{\frac{N-3}{2}}
\exp\Bigl(a\eps-\frac{1}{2}b\eps^2\Bigr) d\eps
\geq
\int\limits_{-\sqrt{N}/2}^{\sqrt{N}/2}\exp\Bigl(
a\eps-\frac{1}{2}L\eps^2
\Bigr)d\eps
\\
&=&
\frac{1}{\sqrt{L}}
\exp\Bigl(\frac{a^2}{2L}\Bigr)
\int\limits_{-\sqrt{LN}/2-a}^{\sqrt{LN}/2-a}
\exp\Bigl(-\frac{x^2}{2}\Bigr)
dx
\geq 
\frac{1}{L}
\int\limits_{-L\sqrt{N}-a}^{L\sqrt{N}-a}
\exp\Bigl(-\frac{x^2}{2}\Bigr)
dx.
\end{eqnarray*}
When $|a|\leq L\sqrt{N}+1,$ this implies that 
$Z_1\geq 1/L.$ Otherwise, say, when $a\geq L\sqrt{N}+1,$
we can use the well known estimates for the Gaussian tail
to write
\begin{eqnarray*}
\int\limits_{-L\sqrt{N}-a}^{L\sqrt{N}-a}
\exp\Bigl(-\frac{x^2}{2}\Bigr)dx
&\geq& 
\frac{1}{L(a-L\sqrt{N})}\exp\Bigl(
-\frac{1}{2}(a-L\sqrt{N})^2
\Bigr)
\\
&&
-L\exp\Bigl(
-\frac{1}{2}(a+L\sqrt{N})^2
\Bigr)
\geq
\frac{1}{L}
\exp(-La^2)
\end{eqnarray*}
which proves (\ref{bound1}). Finally, (\ref{ceps})
and (\ref{bound1}) imply that
$$
\nu_0\bigl(\exp c_0\eps^2\bigr) \leq L\e \exp\bigl(La^2\bigr)
=L\e \exp\bigl(L(z\beta\sqrt{\xi'(q)}+h)^2\bigr)\leq L,
$$
if $\beta$ is small enough, $L\beta^2\xi'(q)<1/2.$
\qed

We are now ready to prove Theorem \ref{epscontrol}.

\textbf{Proof of Theorem \ref{epscontrol}.}
Let us apply (\ref{der}) to $f=\eps^{2k}$ for integer $k\geq 1.$
Since factors $a_l$ and $a_{l,l'}$ are second degree polynomials
in the last coordinates and $|\hR_{l,l'}-q|\leq L,$ 
$|\hR_l - r|\leq L$ we can bound the derivative by
$$
\Bigl| \nu_t'(\eps^{2k}) \Bigr| \leq 
L(\beta^2+h)\nu_t\bigl((1+\eps_1^2+\eps_2^2)\eps_1^{2k}\bigr)
+\nu_t\bigl(\eps^{2k}|\R|\bigr)
\leq L(\beta^2+h)\nu_t\bigl((1+\eps^2)\eps^{2k}\bigr) + \nu_t\bigl(\eps^{2k}|\R|\bigr).
$$
Since $\eps_l^2\leq N,$ for a polynomial $p(\eps_1,\eps_2,\eps_3)$ of the
fourth degree we have
$$
\frac{1}{N}p(\eps_1,\eps_2,\eps_3)\leq\frac{L}{N}\sum_{l\leq 3}(1+\eps_l^4)
\leq L\Bigl(1+\sum_{l\leq 3}\eps_l^2\Bigr).
$$
Therefore,
$$
|\R|\leq L(\beta^2+h)L\Bigl(1+\sum_{l\leq 3}\eps_l^2\Bigr)
$$
and
$$
\Bigl| \nu_t'(\eps^{2k}) \Bigr| 
\leq L(\beta^2+h)\Bigl(\nu_t\bigl((1+\eps^2)\eps^{2k}\bigr)\Bigr).
$$
Using this, we can write
\begin{eqnarray*}
\nu_t'\bigl(\exp c\eps^2\bigr)=\sum_{k\geq 1}\frac{c^k}{k!}\nu_t'(\eps^{2k})
&\leq& 
L(\beta^2+h)\sum_{k\geq 1}\frac{c^k}{k!}\nu_t\bigl((1+\eps^2)\eps^{2k}\bigr)
\\
&\leq& 
L(\beta^2+h)\nu_t\bigl((1+\eps^2)\exp c\eps^2\bigr).
\end{eqnarray*}
If we take $c_0<1$ and let $c(t)=\bigl(c_0-L(\beta^2+h)t\bigr)$ then 
\begin{eqnarray*}
\nu_t'(\exp c(t)\eps^2)
&\leq&
L(\beta^2+h)\nu_t\bigl((1+\eps^2)\exp c(t)\eps^2\bigr)
\\
&-&
L(\beta^2+h)\nu_t\bigl(\eps^2\exp c(t)\eps^2\bigr)
=
L(\beta^2+h)\nu_t\bigl(\exp c(t)\eps^2\bigr).
\end{eqnarray*}
Integrating this over $t$ yields
\begin{equation}
\nu_t(\exp c(t)\eps^2)\leq \exp(L(\beta^2)t)\nu_0(\exp c_0\eps^2)\leq L
\label{epst}
\end{equation}
for small enough $\beta,$ by Lemma \ref{epszero}.
If $\beta^2+h$ is small enough then $c(t)>c_0/2$ and this
finishes the proof of Theorem \ref{epscontrol}.
\qed

\textbf{Proof of Theorem \ref{zeromoments}}
Let us denote 
$$
f(\eps)=\Bigl(1-\frac{\eps^2}{N}\Bigr)^{\frac{N-3}{2}}
\exp\Bigl(a\eps-\frac{1}{2}b\eps^2\Bigr).
$$
Then, using (\ref{decop}) as in (\ref{ceps}), we can write

\begin{eqnarray*}
Z_1 \bigl\la \eps^k\bigr\ra_0
&=&
\int\limits_{-\sqrt{N}}^{\sqrt{N}}
\eps^k f(\eps) d\eps
=
-\frac{1}{b}
\int\limits_{-\sqrt{N}}^{\sqrt{N}}
\eps^{k-1} \Bigl(1-\frac{\eps^2}{N}\Bigr)^{\frac{N-3}{2}}
\exp(a\eps)
d\exp\Bigl(-\frac{1}{2}b\eps^2\Bigr) 
\\
&=&
\frac{1}{b}
\int\limits_{-\sqrt{N}}^{\sqrt{N}}
\Bigl(
(k-1)\eps^{k-2}+a\eps^{k-1} 
\Bigr)
f(\eps) d\eps
-\frac{1}{b}\frac{N-3}{N}
\int\limits_{-\sqrt{N}}^{\sqrt{N}}
\eps^k \Bigl(1-\frac{\eps^2}{N}\Bigr)^{-1}
f(\eps) d\eps
\end{eqnarray*}
by integration by parts.
Moving the last integral to the left hand side of the equation,
\begin{eqnarray}
&&
\int\limits_{-\sqrt{N}}^{\sqrt{N}}
\Bigl(
1+\frac{1}{b}\frac{N-3}{N-\eps^2}
\Bigr)
\eps^k f(\eps) d\eps
=
\frac{1}{b}
\int\limits_{-\sqrt{N}}^{\sqrt{N}}
\Bigl(
(k-1)\eps^{k-2}+a\eps^{k-1} 
\Bigr)
f(\eps) d\eps.
\label{reaba}
\end{eqnarray}
If we rewrite
$$
1+\frac{1}{b}\frac{N-3}{N-\eps^2}
=\frac{b+1}{b}\Bigl(
1+\frac{\eps^2-3}{(b+1)(N-\eps^2)}
\Bigr)
$$
then (\ref{reaba}) implies
\begin{eqnarray*}
\int\limits_{-\sqrt{N}}^{\sqrt{N}}
\eps^k f(\eps) d\eps
&=&
\frac{1}{b+1}
\int\limits_{-\sqrt{N}}^{\sqrt{N}}
\Bigl(
(k-1)\eps^{k-2}+a\eps^{k-1} 
\Bigr)
f(\eps) d\eps
\\
&&
+\frac{1}{N(b+1)}
\int\limits_{\sqrt{N}}^{\sqrt{N}}
\eps^k(3-\eps^2)
\Bigl(1-\frac{\eps^2}{N}\Bigr)^{-1}f(\eps)d\eps.
\end{eqnarray*}
Dividing both sides by $Z_1$ gives
\begin{equation}
S_k=\frac{a}{b+1}S_{k-1}+\frac{k-1}{b+1}S_{k-2}+r_k
\label{Ses}
\end{equation}
where we denoted $S_k=\la\eps^k\ra_0$ and where
$$
r_k=\frac{1}{N(b+1)}\Bigl\la\eps^k(3-\eps^2)
\Bigl(1-\frac{\eps^2}{N}\Bigr)^{-1} \Bigr\ra_0. 
$$
Comparing (\ref{Ses}) with (\ref{gammas}), it should be obvious that 
$S_k=\gamma_k +\hat{r}_k,$ 
where $\hat{r}_k$ is a polynomial in $a$ and $(r_l)_{l\leq k}$ where
each term contains a least one factor $r_l.$ Therefore,
$$
S_{k_1}\ldots S_{k_n}=\gamma_{k_1}\ldots\gamma_{k_n}+r
$$
where $r$ is a polynomial in $a$ and $(r_l)_{l\leq k_0}$
for $k_0=\max(k_1,\ldots,k_n)$ and each term contains at least one
factor $r_l.$ Therefore, each term in $r$ will have at least one
factor $1/N$ and if we can show that for any $k, m>0$
\begin{equation}
\e\Bigl\la
\Bigl(
\eps^{k}(3-\eps^2)
\Bigl(1-\frac{\eps^2}{N}\Bigr)^{-1}
\Bigr)^m
\Bigr\ra_0\leq L
\label{aba}
\end{equation}
then, by H\"older's inequality, $\e|r|\leq L/N$ and this finishes
the proof of Theorem \ref{zeromoments}. To prove (\ref{aba}),
we write that for any polynomial $p(\eps),$ by (\ref{decop}),
$$
\e\Bigl\la
p(\eps)
\Bigl(1-\frac{\eps^2}{N}\Bigr)^{-m}
\Bigr\ra_0
=
\e \frac{1}{Z_1}\int\limits_{-\sqrt{N}}^{\sqrt{N}}
p(\eps) 
\Bigl(1-\frac{\eps^2}{N}\Bigr)^{\frac{N-3-2m}{2}}
\exp\Bigl(a\eps-\frac{1}{2}b\eps^2\Bigr)d\eps.
$$
Repeating the argument of Lemma \ref{epszero} one can show
that for small enough $\beta>0$ the right hand side is bounded
by some $L>0$ which proves (\ref{aba}). 
\qed

\section{Control of the overlap and magnetization.}\label{SecR}

We finally turn to the proof of Theorem \ref{Big}.
We will start with the following result.
Given a set $A\subseteq S_{N-1}^n,$ let us denote
$$
I_A=I\bigl((\hS^1,\ldots,\hS^n)\in A\bigr).
$$
Then the following Lemma holds.

\begin{lemma}\label{tint}
If $A\subseteq S_{N-1}^n$ is symmetric with respect to permutations
of the coordinates, then for small enough $\beta$ and $h,$
\begin{equation}
\Bigl|
\frac{1}{N}\e\log\bigl\la I_A \bigr\ra_t
-\frac{1}{N}\e\log\bigl\la I_A \bigr\ra_0
\Bigr|\leq \frac{L}{N}
\label{phit}
\end{equation}
\end{lemma}
We will apply (\ref{phit}) to the sets of the type
\begin{equation}
\Bigl\{\hS^1 : |\hR_1 -r|\geq x \Bigr\} 
\,\,\,\mbox{ or }\,\,\,
\Bigl\{(\hS^1,\hS^2) : |\hR_{1,2} -q|\geq x \Bigr\}
\label{Atype}
\end{equation}
and Lemma \ref{tint} states that their Gibbs' measure does not
change much along the interpolation (\ref{Hamint}).

\textbf{Proof of Lemma \ref{tint}.}
For a set $A\subseteq S_{N-1}^n,$ let us consider
$$
\phi_A(t)=\frac{1}{N}\e \log 
\int_{S_{N}^n}I_{A}\exp\sum_{l\leq n}H_t(\vsi^l)d\lambda_N^n.
$$
Then 
$$
\frac{1}{N}\e\log\bigl\la I_A \bigr\ra_t
=\phi_A(t)-\phi_{S_{N-1}^n}(t)
$$
and Lemma \ref{tint} follows from the following.

\begin{lemma}\label{intcontrol}
For small enough $\beta$ and $h$ we have
\begin{equation}
|\phi_A'(t)|\leq \frac{L}{N}.
\label{FAder}
\end{equation}
\end{lemma}
\textbf{Proof.}
Given  a function $f=f(\vsi^1,\ldots,\vsi^n),$ we define
\begin{equation}
\la f\ra_{t,A}= \frac{\la fI_A\ra_t}{\la I_A\ra_t} =
\int_{S_N^n}I_A f \exp\sum_{l\leq n}H_t(\vsi^l)d\lambda_N^n
\Bigr/\int_{S_N^n}I_A \exp\sum_{l\leq n}H_t(\vsi^l)d\lambda_N^n.
\label{GA}
\end{equation}
Then
$$
N\phi_A'(t)=\e\Bigl\la \sum_{l\leq n} \frac{\partial H_t(\vsi^l)}{\partial t}\Bigr\ra_{t,A}.
$$
If we denote
\begin{eqnarray*}
&&
S_{l,l'} =N \xi(R_{l,l'}) -(N-1)\xi(\hR_{l,l'})
-\eps_l\eps_{l'} \xi'(q)
\end{eqnarray*}
then integration by parts as in Theorem \ref{Thder} gives,
\begin{eqnarray}
N\phi_A'(t)
&=&
\sum_{l\leq n}\e\Bigl\la
h\sum_{i\leq N-1}\hs_i^l\Bigl(
\sqrt{\frac{N-\eps_l^2}{N-1}}-1
\Bigr)+\frac{1}{2}\eps_l^2 b
\Bigr\ra_{t,A}
\nonumber
\\
&+&
\frac{\beta^2}{2}\sum_{l,l'\leq n}\e\la S_{l,l'}\ra_{t,A}
-\frac{\beta^2}{2}\sum_{l\leq n}\sum_{l'=n+1}^{2n}
\e\la S_{l,l'}\ra_{t,A}.
\label{derA}
\end{eqnarray}
The Gibbs average in the last term is defined on two copies
$(\vsi^1,\ldots,\vsi^n)$ and $(\vsi^{n+1},\ldots,\vsi^{2n}).$
Since
$$
\Bigl|
(N-1)\Bigl(
\sqrt{\frac{N-\eps_l^2}{N-1}}-1
\Bigr)
\Bigr|
\leq
L(1+\eps_l^2)
$$
and $|S_{l,l'}|\leq L(1+\eps_l^2+\eps_{l'}^2),$ (\ref{derA}) implies that
\begin{equation}
\bigl| N\phi_A'(t)\bigr|\leq L(1+\sum_{l\leq n}\e\la \eps_l^2\ra_{t,A})
\leq L(1+\e\la \eps_1^2\ra_{t,A}),
\label{aaa}
\end{equation}
where in the last inequality we used the fact that $\e\la \eps_l^2\ra_{t,A}$
does not depend on $l$ due to the symmetry of $A.$
One can now repeat the proof of Theorem \ref{epscontrol} to obtain
the analogue of (\ref{epst}):
$$
\e\bigl\la \exp c(t) \eps_1^2\bigr\ra_{t,A} \leq \exp\bigl(L(\beta^2+h)t\bigr)
\e\bigl\la \exp c_0\eps_1^2\bigr\ra_{0,A},
$$
where $c(t)=c_0-L(\beta^2+h)t>c_0/2$ for small enough $\beta,h.$
Using (\ref{GA}) and (\ref{decop}), we can write
$$
\e\bigl\la \exp c_0\eps_1^2\bigr\ra_{0,A}
=
\e\frac{\bigl\la I_A \exp c_0\eps_1^2\bigr\ra_{0}}
{\bigl\la I_A\bigr\ra_{0}}
=
\e\frac{\bigl\la I_A\bigr\ra_0 \bigl\la\exp c_0\eps_1^2\bigr\ra_{0}}
{\bigl\la I_A\bigr\ra_{0}}
=
\e \bigl\la\exp c_0\eps_1^2\bigr\ra_{0} \leq L
$$
for $c_0<1$ and small enough $\beta,$ by Lemma \ref{epszero}. 
Hence, $\e\bigl\la \eps_1^2\bigr\ra_{t,A}\leq L$ and
(\ref{aaa}) finishes the proof of Lemma \ref{intcontrol}.
\qed

To apply Lemma \ref{tint} to the sets of the type (\ref{Atype}),
we need to control $N^{-1}\e\log\bigl\la I_A\bigr\ra_{0}.$ Let us notice
that $\bigl\la I_A\bigr\ra_0$ for the sets in (\ref{Atype}) is defined
exactly in the same way as $\bigl\la I_A\bigr\ra$ (i.e. for $t=1$)
for the sets of the type
\begin{equation}
\Bigl\{\vsi^1 : |R_1 -r|\geq x \Bigr\} 
\,\,\,\mbox{ or }\,\,\,
\Bigl\{(\vsi^1,\vsi^2) : |R_{1,2} -q|\geq x \Bigr\}
\label{Atype1}
\end{equation}
only for $N-1$ instead of $N.$ Therefore, for simplicity of notations,
we will show how to control $N^{-1}\e\log\bigl\la I_A\bigr\ra$ for $A$
in (\ref{Atype1}) and then apply it to (\ref{Atype}).

For $\bq\in[0,1]$ consider a Hamiltonian
\begin{equation}
h_t(\vsi)=\sqrt{t}H_N(\vsi)+\sum_{i\leq N}\sigma_i \Bigl(
\sqrt{1-t}z_i\beta\sqrt{\xi'(\bq)}+h
\Bigr).
\label{ht}
\end{equation}
Let $\bla\cdot\bra_t$ define the Gibbs average with respect to the Hamiltonian
(\ref{ht}). Let us define $\bq$ as any solution of the
equation
\begin{equation}
\bq=\e\bla R_{1,2}\bra_0
\label{qbar}
\end{equation}
where the right hand side depends on $\bq$ through (\ref{ht}). We will show
that there exists a solution close to $q.$ Given $\bq$ that satisfies (\ref{qbar})
we define
\begin{equation}
\bar{r}=\e\bla R_1\bra_0.
\label{rbar}
\end{equation}

\begin{lemma}\label{qrbar} For small enough $\beta, h$ 
there exists a solution of (\ref{qbar}) such that
$$
|\bq -q|\leq \frac{L\log^2 N}{N},\,\,\,
|\bar{r} -r|\leq \frac{L\log^2 N}{N}.
$$
\end{lemma}

We will also prove  the following.

\begin{lemma}\label{RandR}
For small enough $\beta$ we can find $\alpha>0$ such that
for $\bar{q}, \bar{r}$ as in Lemma \ref{qrbar}, 
$$
\e\bigl\la
\exp N\alpha(R_{1,2}-\bq)^2
\bigr\ra \leq L 
\,\,\,\mbox{ and }\,\,\,
\e\bigl\la
\exp N\alpha(R_{1}-\bar{r})^2
\bigr\ra \leq L. 
$$
\end{lemma}

Before we prove Lemmas \ref{qrbar} and \ref{RandR}, let us first show
how they together with Lemma \ref{tint} imply Theorem \ref{Big}.

\textbf{Proof of Theorem \ref{Big}.}
Lemma \ref{RandR} implies that
\begin{eqnarray*}
\e\log\Bigl\la
I\bigl(|R_{1,2}-\bq|\geq x\bigr)
\Bigr\ra 
&\leq &
\e\log\Bigl\la
\exp N\alpha (R_{1,2}-\bq)^2
\Bigr\ra -N\alpha x^2
\\
&\leq&
\log\e\Bigl\la
\exp N\alpha (R_{1,2}-\bq)^2
\Bigr\ra -N\alpha x^2
\leq L-N\alpha x^2.
\end{eqnarray*}
Using this for $N-1$ instead of $N$ yields
$$
\frac{1}{N}\e\log\Bigl\la
I\bigl(|\hR_{1,2}-\bq|\geq x\bigr)
\Bigr\ra_0 
\leq \frac{L}{N}-\alpha x^2
$$
and by Lemma \ref{tint}
$$
\frac{1}{N}\e\log\Bigl\la
I\bigl(|\hR_{1,2}-\bq|\geq x\bigr)
\Bigr\ra_t
\leq \frac{L}{N}-\alpha x^2.
$$
For $x=L\bigl(\log N/N\bigr)^{1/4}$ we get
$$
\frac{1}{N}\e\log\Bigl\la
I\Bigl(|\hR_{1,2}-\bq|\geq L\Bigl(
\frac{\log N}{N}\Bigr)^{1/4}\Bigr)
\Bigr\ra_t
\leq -L\Bigl(\frac{\log N}{N}\Bigr)^{1/2}=:\delta.
$$
Gaussian concentration of measure (as in Corollary 2.2.5 in \cite{SG})
implies that
$$
\frac{1}{N}\log\Bigl\la
I\Bigl(|\hR_{1,2}-\bq|\geq L\Bigl(
\frac{\log N}{N}\Bigr)^{1/4}\Bigr)
\Bigr\ra_t
\leq -L\Bigl(\frac{\log N}{N}\Bigr)^{1/2}.
$$
with probability at least $1-L\exp(-N\delta^2/L)\geq 1-LN^{-K}$
for any $K>0,$ by choosing $L$ in the definition of $x$ sufficiently
large. Therefore, with probability at least $1-LN^{-K},$
$$
\Bigl\la
I\Bigl(|\hR_{1,2}-\bq|\geq L\Bigl(
\frac{\log N}{N}\Bigr)^{1/4}\Bigr)
\Bigr\ra_t
\leq \exp\Bigl(-L\bigl(N\log N\bigr)^{1/2}\Bigr)
\leq LN^{-K}
$$
and, thus, 
$$
\e\Bigl\la
I\Bigl(|\hR_{1,2}-\bq|\geq L\Bigl(
\frac{\log N}{N}\Bigr)^{1/4}\Bigr)
\Bigr\ra_t
\leq LN^{-K}
$$
Lemma \ref{qrbar} implies
$$
\e\Bigl\la
I\Bigl(|\hR_{1,2}-q|\geq L\Bigl(
\frac{\log N}{N}\Bigr)^{1/4}\Bigr)
\Bigr\ra_t
\leq LN^{-K}
$$
and this proves the first part of Theorem \ref{Big}. 
The second part is proved similarly.
\qed

\textbf{Proof of Lemma \ref{qrbar}.}
If we denote
$$
\vec{v}=\Bigl(z_1\beta\sqrt{\xi'(\bq)}+h,\ldots,z_N \beta\sqrt{\xi'(\bq)}+h\Bigr)
$$
Then
$$
\bla R_{1,2}\bra_0=\frac{1}{Z^2}\int\limits_{S_N^2}\frac{1}{N}(\vsi^1,\vsi^2)
\exp\bigl((\vsi^1,\vec{v})+(\vsi^2,\vec{v})\bigr)d\lambda_N(\vsi^1)d\lambda_N(\vsi^2),
$$
where $Z=\int_{S_N}\exp((\vsi,\vec{v}))d\lambda_N(\vsi).$
If $O$ is an orthogonal transformation such that 
$O\vec{v}=(0,\ldots,0,|\vec{v}|)$ then making a change of variables
$\vsi^l\to O^{-1}\vsi^l$ we get
$$
\bla R_{1,2}\bra_0 =\frac{1}{Z^2}\int\limits_{S_N^2}\frac{1}{N}(\vsi^1,\vsi^2)
\exp\bigl(\eps_1|\vec{v}|+\eps_2|\vec{v}|\bigr)d\lambda_N(\vsi^1)d\lambda_N(\vsi^2),
$$
and $Z=\int_{S_N}\exp(\eps|\vec{v}|)d\lambda_N(\vsi).$
By (\ref{ars}) 
$$
\frac{1}{N}(\vsi^1,\vsi^2)=R_{1,2}=
\sqrt{\Bigl(1-\frac{\eps_1^2}{N}\Bigr)\Bigl(1-\frac{\eps_2^2}{N}\Bigr)} \hR_{1,2}
+\frac{1}{N}\eps_1\eps_2
$$
and by (\ref{Fub2})
\begin{eqnarray*}
&&
\int\limits_{S_N^2}
\sqrt{\Bigl(1-\frac{\eps_1^2}{N}\Bigr)\Bigl(1-\frac{\eps_2^2}{N}\Bigr)} \hR_{1,2}
\exp\bigl(\eps_1|\vec{v}|+\eps_2|\vec{v}|\bigr)d\lambda_N(\vsi^1)d\lambda_N(\vsi^2)
\\
&&
=
a_N^2 \Bigl(
\int\limits_{-\sqrt{N}}^{\sqrt{N}}
\Bigl(1-\frac{\eps^2}{N}\Bigr)^{\frac{N-2}{2}} \exp(\eps|\vec{v}|)d\eps
\Bigr)^2
\int\limits_{S_{N-1}^2}\hR_{1,2}d\lambda_{N-1}(\hS^1)d\lambda_{N-1}(\hS^2)
=0
\end{eqnarray*}
since the last integral is equal to zero by symmetry. Therefore,
\begin{equation}
\bla R_{1,2}\bra_0
=
\frac{1}{N}\bla \eps_1\eps_2 \bra_0
=
\bla N^{-1/2}\eps {\bra_0}^2
\label{Reps}
\end{equation}
and using (\ref{Fub2}) again
$$
\bla N^{-1/2}\eps\bra_0
=
\int\limits_{-\sqrt{N}}^{\sqrt{N}}
\frac{\eps}{\sqrt{N}}
\Bigl(1-\frac{\eps^2}{N}\Bigr)^{\frac{N-3}{2}} \exp(\eps|\vec{v}|)d\eps
\Bigr/
\int\limits_{-\sqrt{N}}^{\sqrt{N}}
\Bigl(1-\frac{\eps^2}{N}\Bigr)^{\frac{N-3}{2}} \exp(\eps|\vec{v}|)d\eps.
$$
By making a change of variable $\eps=\sqrt{N}x$ we can rewrite the right hand side as
\begin{equation}
\bla N^{-1/2}\eps\bra_0
=
\int\limits_{-1}^{1}x \exp N\varphi(x) dx
\Bigr/
\int\limits_{-1}^{1} \exp N\varphi(x) dx
\label{upthere}
\end{equation}
where
\begin{equation}
\varphi(x)=cx+\frac{N-3}{2N}\log(1-x^2)
\label{fi}
\end{equation}
and
\begin{equation}
c=N^{-1/2}|\vec{v}|=\Bigl(
\frac{1}{N}\sum_{i\leq N} \bigl( z_i\beta\sqrt{\xi(\bq)}+h\bigr)^2
\Bigr)^{1/2}.
\label{cee}
\end{equation}
Let $x_0$ denotes the point where $\varphi(x)$ achieves its maximum
which satisfies
\begin{equation}
\varphi'(x_0)=0\Longrightarrow
c=\frac{N-3}{N}\frac{x_0}{1-x_0^2}.
\label{cx}
\end{equation}
Since $|\eps|/\sqrt{N}\leq 1$ and $|x_0|\leq 1,$
\begin{equation}
\Bigl|
\e\bla N^{-1/2}\eps{\bra_0}^2
-\e x_0^2
\Bigr| \leq 2\e
\Bigl|
\int\limits_{-1}^{1}(x-x_0) \exp N\varphi(x) dx
\Bigr|
\Bigr/
\int\limits_{-1}^{1} \exp N\varphi(x) dx.
\label{dif}
\end{equation}
For $c$ in (\ref{cee}) and $c'>2h^2,$
\begin{eqnarray}
\p\bigl(c\geq c'\bigr)
&=&
\p\Bigl(
\sum_{i\leq N} \bigl(z_i\beta\sqrt{\xi'(\bq)}+h\bigr)^2
\geq N{c'}^2
\Bigr)
\label{cc}
\\
&\leq&
\p\Bigl(
2\beta^2\xi'(\bq)\sum_{i\leq N} z_i^2
\geq N\bigl({c'}^2 -2h^2\bigr)
\Bigr)
=
\p\Bigl(\sum_{i\leq N}z_i^2 \geq Nc''\Bigr)
\leq
\exp(-LN),
\nonumber
\end{eqnarray}
where $L$ can be made arbitrarily large by increasing $c'.$

Let us now assume that the event $\{c\leq c'\}$ occurs. Then (\ref{cx})
implies that $|x_0|\leq 1-\delta$ for some $\delta>0$ that depends on
$c'$ only. Let us define
$$
\Omega=\Bigl\{
x\in[-1,1] : |x-x_0|\leq \omega = \sqrt{\frac{L\log N}{N}}
\Bigr\}
$$
for $L$ large enough
and write $\int_{-1}^{1}\exp N\varphi(x) dx = \mbox{I} + \mbox{II},$ where
$$
\mbox{I}=\int\limits_{\Omega} \exp N\varphi(x)dx
\,\,\,\mbox{ and }\,\,\,
\mbox{II}=\int\limits_{\Omega^c} \exp N\varphi(x)dx.
$$
We have
\begin{equation}
\varphi''(x)=-\frac{1+x^2}{(1-x^2)^2}\leq -1
\label{derr}
\end{equation}
and for $|x|\leq 1-\delta/2,$ clearly, $-L\leq \varphi''(x)$
and $|\varphi'''(x)|\leq L.$
Since $\varphi'(x_0)=0,$ we have $\varphi(x)\geq \varphi(x_0)-L(x-x_0)^2$
for $x\in\Omega$ and, therefore,
\begin{eqnarray}
\mbox{I} 
&\geq& 
\exp N\varphi(x_0)\int\limits_{\Omega}
\exp\Bigl(-LN(x-x_0)^2\Bigr)
\label{uno}
\\
&=&
\exp (N\varphi(x_0))\frac{1}{\sqrt{N}}
\int\limits_{|y|\leq (L\log N)^{1/2}}
\exp\bigl(-Ly^2\bigr)dy
\geq
\frac{1}{L\sqrt{N}}\exp N\varphi(x_0).
\nonumber
\end{eqnarray}
On the other hand, by (\ref{derr}), $\varphi(x)\leq \varphi(x_0)
-(x-x_0)^2/2$ and, thus,
\begin{eqnarray*}
&&
\mbox{II} 
\leq
\exp (N\varphi(x_0))\frac{1}{\sqrt{N}}
\int\limits_{|y|\geq (L\log N)^{1/2}}
\exp\bigl(-Ly^2\bigr)dy
\leq
\frac{L}{N^K}\exp N\varphi(x_0),
\end{eqnarray*}
where $K$ can be made arbitrarily large by a proper choice
of $L$ in the definition of $\Omega.$
The denominator in (\ref{dif}) can be bounded from below by
\begin{equation}
\int\limits_{-1}^{1}\exp N\varphi(x)dx \geq \mbox{I} \geq
\frac{1}{L\sqrt{N}}\exp N\varphi(x_0).
\label{Zlow}
\end{equation}
Next, we write $\int_{-1}^{1}(x-x_0)\exp N\varphi(x) dx = \mbox{III} + \mbox{IV},$ where
$$
\mbox{III}=\int\limits_{\Omega}(x-x_0) \exp N\varphi(x)dx
\,\,\,\mbox{ and }\,\,\,
\mbox{IV}=\int\limits_{\Omega^c} (x-x_0) \exp N\varphi(x)dx.
$$
We control IV by
\begin{equation}
|\mbox{IV}|\leq 2|\mbox{II}|\leq \frac{L}{N^{K}}\exp N\varphi(x_0).
\label{four}
\end{equation}
To control III we use that for $x\in\Omega$
$$
\Bigl|
\varphi(x)-\varphi(x_0)-\frac{1}{2}\varphi''(x_0)(x-x_0)^2
\Bigr|
\leq L\Bigl(\frac{\log N}{N}\Bigr)^{3/2}=:\Delta.
$$
We have
\begin{eqnarray*}
\mbox{III}
&=&
\int\limits_{x_0-\omega}^{x_0}(x-x_0)\exp N\varphi(x)dx
+\int\limits_{x_0}^{x_0+\omega}(x-x_0)\exp N\varphi(x)dx
\\
&\leq& 
\int\limits_{x_0-\omega}^{x_0}(x-x_0)\exp N\Bigl(
\varphi(x_0)+\frac{1}{2}\varphi''(x_0)(x-x_0)^2-\Delta
\Bigr)dx
\\
&+&
\int\limits_{x_0}^{x_0+\omega}(x-x_0)\exp N\Bigl(
\varphi(x_0)+\frac{1}{2}\varphi''(x_0)(x-x_0)^2+\Delta
\Bigr)dx
\\
&=&
(e^{N\Delta}-e^{-N\Delta})
\int\limits_{x_0}^{x_0+\omega}(x-x_0)\exp N\Bigl(
\varphi(x_0)+\frac{1}{2}\varphi''(x_0)(x-x_0)^2
\Bigr)dx
\\
&\leq&
LN\Delta\omega\exp N\varphi(x_0)
\int\limits_{x_0}^{x_0+\omega}\exp \Bigl(
-\frac{1}{2}N(x-x_0)^2
\Bigr)dx
\\
&\leq&
LN^{1/2}\Delta\omega\exp N\varphi(x_0)
\leq 
\frac{L\log^2 N}{N^{3/2}}\exp N\varphi(x_0).
\end{eqnarray*}
The lower bound can be carried out similarly and, thus,
$$
|\mbox{III}|\leq \frac{L\log^2 N}{N^{3/2}}\exp N\varphi(x_0).
$$
Combining this with (\ref{dif}), (\ref{cc}), (\ref {uno}) and (\ref{four}) proves
\begin{eqnarray}
\Bigl|
\e\bla N^{-1/2}\eps{\bra_0}^2
-\e x_0^2
\Bigr| 
&\leq &
\exp(-LN) + \frac{L}{N^K}
\nonumber
\\
&+&
\e 
\frac{L\log^2 N \exp N\varphi(x_0)}{N^{3/2}}
\Bigr/ 
\frac{\exp N\varphi(x_0)}{L\sqrt{N}}\leq
\frac{L\log^2 N}{N}.
\nonumber
\end{eqnarray}
By (\ref{Reps}), we proved that
\begin{equation}
\Bigl|
\e\bla R_{1,2}\bra_0
-\e x_0^2
\Bigr| \leq 
\frac{L\log^2 N}{N}.
\label{dif2}
\end{equation}
If we denote
$$
c_N=\frac{N}{N-3}c = \frac{N}{N-3}\frac{|\vec{v}|}{\sqrt{N}}
$$
then solving (\ref{cx}) for $x_0$ gives
\begin{equation}
x_0=\frac{2c_N}{1+\sqrt{1+4c_N^2}} 
\,\,\,\mbox{ and }\,\,\,
x_0^2=1-\frac{2}{1+\sqrt{1+2c_N^2}}.
\label{xx}
\end{equation}
It is easy to check that the first two derivatives of
$y(x)=1/(1+\sqrt{1+4x})$ are bounded by an absolute constant
for $x\geq 0$ and, therefore,
$$
|y(c_N^2)-y(\e c_N^2)-y'(\e c_N^2)(c_N^2-\e c_N^2)| \leq L(c_N^2-\e c_N^2)^2.
$$
Taking expectations proves that
\begin{equation}
\Bigl|
\e x_0^2 - \Bigl(
1-\frac{2}{1+\sqrt{1+4\e c_N^2}}
\Bigr)
\Bigr| \leq L\e(c_N^2-\e c_N^2)^2\leq \frac{L}{N}
\label{tam}
\end{equation}
since
$$
c_N^2=\Bigl(\frac{N}{N-3}\Bigr)^2\frac{1}{N}\sum_{i\leq N}
\bigl(z_i \beta\sqrt{\xi'(\bq)}+h\bigr)^2.
$$
If we denote 
$$
\delta=\e\bla R_{1,2}\bra_0
-\Bigl(
1-\frac{2}{1+\sqrt{1+4\e c_N^2}}
\Bigr)
$$
then (\ref{dif2}) and (\ref{tam}) imply that
$|\delta|\leq L \log^2 N/N.$ By (\ref{qbar}),
$\e\bla R_{1,2}\bra_0 =\bq$ and, therefore,
$$
\bq -\delta = 1-\frac{2}{1+\sqrt{1+4\e c_N^2}}
$$
or, equivalently,
$$
\e c_N^2 = \frac{\bq -\delta}{(1-\bq +\delta)^2}
= 
\Bigl(\frac{N}{N-3}\Bigr)^2(\beta^2\xi'(\bq)+h^2).
$$
Comparing with (\ref{critq}), it is now a simple exercise to show
that
$$
|\bq - q|\leq \frac{L\log^2 N}{N}
$$
and this proves the first part of Lemma \ref{qrbar}. The computation
of $\bar{r}$ is slightly different. 
If $\vec{1}=(1,\ldots,1)\in\Reals^N$ then
$$
\bla R_{1}\bra_0=\frac{1}{Z}\int\limits_{S_N}\frac{1}{N}(\vsi,\vec{1})
\exp (\vsi,\vec{v}) d\lambda_N(\vsi)
=
\frac{1}{Z}\int\limits_{S_N}\frac{1}{N}(O^T\vsi,\vec{1})
\exp \eps|v| d\lambda_N(\vsi),
$$
where $O$ is the orthogonal transformation as above.
Note that the last row of $O$ is $\vec{v}/|\vec{v}|.$
Next, we use (\ref{Fub1}) to write
$
\int_{S_N}(O^T\vsi,\vec{1})
\exp \eps|v| d\lambda_N(\vsi)
$
as
\begin{eqnarray*}
&&
a_N \int\limits_{-\sqrt{N}}^{\sqrt{N}}d\eps 
\exp\eps|\vec{v}| \Bigl(
1-\frac{\eps^2}{N}
\Bigr)^{\frac{N-3}{2}}
\int\limits_{S_{N-1}}
\Bigl(
O^T\Bigl(\sqrt{\frac{N-\eps^2}{N-1}}\hS,\eps\Bigr),
\vec{1}\Bigr)d\lambda_{N-1}(\hS).
\\
&&
=
a_N \int\limits_{-\sqrt{N}}^{\sqrt{N}}d\eps 
\exp\eps|\vec{v}| \Bigl(
1-\frac{\eps^2}{N}
\Bigr)^{\frac{N-3}{2}}
\int\limits_{S_{N-1}}
\bigl(
O^T\bigl(0,\ldots,0,\eps\bigr),
\vec{1}\bigr)d\lambda_{N-1}(\hS)
\end{eqnarray*}
by symmetry $\hS\to -\hS.$
Since the last column of $O^T$ is $\vec{v}/|\vec{v}|$
$$
\bigl(
O^T\bigl(0,\ldots,0,\eps\bigr),
\vec{1}\bigr)
=
\frac{1}{|\vec{v}|}\eps\sum_{i\leq N} v_i
$$
and, therefore,
$$
\int_{S_N}(O^T\vsi,\vec{1})
\exp \eps|v| d\lambda_N(\vsi)
=a_N 
\frac{1}{|\vec{v}|}\sum_{i\leq N} v_i
\int\limits_{-\sqrt{N}}^{\sqrt{N}}
\eps\exp\eps|\vec{v}| \Bigl(
1-\frac{\eps^2}{N}
\Bigr)^{\frac{N-3}{2}} d\eps.
$$
Similarly
$$
Z = a_N 
\int\limits_{-\sqrt{N}}^{\sqrt{N}}
\exp\eps|\vec{v}| \Bigl(
1-\frac{\eps^2}{N}
\Bigr)^{\frac{N-3}{2}} d\eps
$$
and making the change of variable $\eps=\sqrt{N}x$ we get
$$
\bla R_1 \bra_0 =\frac{1}{\sqrt{N}}
\frac{1}{|\vec{v}|}\sum_{i\leq N} v_i
\int\limits_{-1}^{1}x\exp N\varphi(x) dx
\Bigr/
\int\limits_{-1}^{1}
\exp N\varphi(x)dx.
$$
Repeating the argument leading to (\ref{dif2}) one can now show that
\begin{equation}
\Bigl|
\e\bla R_{1}\bra_0
-\e \frac{1}{\sqrt{N}}
\frac{1}{|\vec{v}|}\sum_{i\leq N} v_i x_0
\Bigr| \leq 
\frac{L\log^2 N}{N}.
\label{difff}
\end{equation}
By (\ref{xx}),
$$
\frac{1}{\sqrt{N}}
\frac{1}{|\vec{v}|}\sum_{i\leq N} v_i x_0
=\frac{1}{N-3}\sum_{i\leq N} v_i \frac{2}{1+\sqrt{1+4c_N^2}}.
$$
Since $c_N^2$ is concentrated near $\e (z_1\beta\sqrt{\xi'(\bq)}+h)^2 =\beta^2 \xi'(\bq) +h^2$
and $\e  v_i = h,$ it is a simple exercise to show that
$$
\Bigl|
\e \frac{1}{N-3}\sum_{i\leq N} v_i \frac{2}{1+\sqrt{1+4c_N^2}}
-\frac{2h}{1+\sqrt{1+4(\beta^2 \xi'(\bq) + h^2)}}
\Bigr| \leq \frac{L}{N}.
$$
Since $|\bq -q|\leq L\log^2 N/N,$ we get
$$
\Bigl|
\e\bla R_{1}\bra_0
-\frac{2h}{1+\sqrt{1+4(\beta^2 \xi'(q) + h^2)}}
\Bigr| \leq \frac{L\log^2 N}{N}
$$
and since by (\ref{critq})
$$
\frac{2h}{1+\sqrt{1+4(\beta^2 \xi'(q) + h^2)}}=h(1-q)=r
$$
we proved that $|\bar{r} - r|\leq L\log^2 N/N.$
This finishes the proof of Lemma \ref{qrbar}.
\qed

\textbf{Proof of Lemma \ref{RandR}.}
We notice that $\bigl\la\cdot \bigr\ra = \bla\cdot\bra_1$
so the proof will proceed by interpolation in (\ref{ht}).
If is easy to show similarly to Theorem \ref{Thder} that
for a function $f=f(\vsi^1,\ldots,\vsi^n),$
\begin{eqnarray*}
\frac{\partial}{\partial t}\e\bla f\bra_t
&=& 
N\beta^2 \sum_{1\leq l<l'\leq n}
\e \bla f\Delta(R_{l,l'})\bra_t 
-N\beta^2n\sum_{l\leq n}
\e \bla f\Delta(R_{l,n+1})\bra_t
\\
&&
+
N\beta^2\frac{n(n+1)}{2}\e \bla f\Delta(R_{n+1,n+2})\bra_t,
\end{eqnarray*}
where
$$
\Delta(R_{l,l'})=\xi(R_{l,l'})-R_{l,l'}\xi'(\bq)+\theta(\bq)
$$
and $\theta(x)=x\xi'(x)-\xi(x).$ Since $\xi$ is three times continuously
differentiable we have
$$
|\Delta(R_{l,l'})| \leq L(R_{l,l'}-\bq)^2.
$$
For $n=2$ and for any $k\geq 1$ this implies, by H\"older's inequality,
$$
\frac{\partial}{\partial t}\e\bla (R_{1,2}-\bq)^{2k}\bra_t
\leq L N\beta^2 \e\bla (R_{1,2}-\bq)^{2k+2}\bra_t.
$$ 
Therefore, 
\begin{eqnarray*}
\frac{\partial}{\partial t}\e\bla \exp N \alpha(R_{1,2}-\bq)^{2}\bra_t
&\leq&
\sum_{k\geq 1} LN \beta^2 \frac{N^k \alpha^k}{k!}
\e\bla (R_{1,2}-\bq)^{2k+2}\bra_t
\\
&\leq& 
L N\beta^2 \e\bla (R_{1,2}-\bq)^2 \exp N\alpha(R_{1,2}-\bq)^{2}\bra_t.
\end{eqnarray*}
For $\alpha(t)=\alpha - L \beta^2 t$ this implies
$$
\frac{\partial}{\partial t}\e\bla \exp N \alpha(t)(R_{1,2}-\bq)^{2}\bra_t
\leq 0
$$
and, therefore,
$$
\e\bla \exp N \alpha(t)(R_{1,2}-\bq)^{2}\bra_t
\leq
\e\bla \exp N \alpha(R_{1,2}-\bq)^{2}\bra_0.
$$
Next, since
\begin{eqnarray*}
\e\bla (R_1 - \bar{r})^{2k}(R_{l,l'}-\bq)^{2}\bra_t
&\leq&
\Bigl(\e\bla (R_{1}-\bar{r})^{2k+2}\bra_t\Bigr)^{\frac{2k}{2k+2}}
\Bigl(\e\bla (R_{1,2}-\bq)^{2k+2}\bra_t\Bigr)^{\frac{2}{2k+2}}
\\
&\leq&
\frac{k}{k+1}\e\bla (R_{1}-\bar{r})^{2k+2}\bra_t
+\frac{1}{k+1}\e\bla (R_{1,2}-\bq)^{2k+2}\bra_t
\end{eqnarray*}
we can bound the derivative of
$\e\bla \exp N \alpha(R_{1}-\bar{r})^{2}\bra_t$
by
\begin{eqnarray*}
&&
LN\beta^2
\sum_{k\geq 1} \frac{N^k \alpha^k}{k!}
\Bigl(
\frac{k}{k+1}\e\bla (R_{1}-\bar{r})^{2k+2}\bra_t
+\frac{1}{k+1}\e\bla (R_{1,2}-\bq)^{2k+2}\bra_t
\Bigr)
\\
&&
\leq
LN\beta^2
\e\bla (R_1-\bar{r})^2\exp N \alpha(R_{1}-\bar{r})^{2}\bra_t
+\frac{L\beta^2}{\alpha} 
\e\bla \exp N \alpha(R_{1,2}-\bq)^{2}\bra_t.
\end{eqnarray*}
For $\alpha(t)=\alpha -L\beta^2 t$ this implies that
\begin{eqnarray*}
\frac{\partial}{\partial t}\e\bla \exp N \alpha(t)(R_{1}-\bar{r})^{2}\bra_t
&\leq&
\frac{L\beta^2}{\alpha(t)} 
\e\bla \exp N \alpha(t)(R_{1,2}-\bq)^{2}\bra_t
\\
&\leq&
\frac{L\beta^2}{\alpha-L\beta^2} 
\e\bla \exp N \alpha(R_{1,2}-\bq)^{2}\bra_0
\end{eqnarray*}
and, thus,
$$
\e\bla \exp N \alpha(1)(R_{1}-\bar{r})^{2}\bra_1
\leq
\e\bla \exp N \alpha(R_{1}-\bar{r})^{2}\bra_0
+
\frac{L\beta^2}{\alpha-L\beta^2} 
\e\bla \exp N \alpha(R_{1,2}-\bq)^{2}\bra_0.
$$
To finish the proof of Lemma \ref{RandR} it remains to show that
for small enough $\alpha,$
$$
\e\bla \exp N \alpha(R_{1,2}-\bq)^{2}\bra_0 \leq L
\,\,\,\mbox{ and }\,\,\,
\e\bla \exp N \alpha(R_{1}-\bar{r})^{2}\bra_0 \leq L.
$$
By (\ref{qbar}) and Jensen's inequality
\begin{eqnarray*}
\e\bla \exp N \alpha(R_{1,2}-\bq)^{2}\bra_0 
&\leq&
\e\bla \exp N \alpha(R_{1,2}-R_{3,4})^{2}\bra_0. 
\\
&=&
\e\frac{1}{Z^4}
\int\limits_{S_{N}^4}
\exp\Bigl(
N\alpha(R_{1,2}-R_{3,4})^2+\sum_{l\leq 4}(\vsi^l,\vec{v})
\Bigr)d\lambda_N^4,
\end{eqnarray*}
as in the beginning of Lemma \ref{qrbar}.
For $\vec{v}$ and $O$ defined in Lemma \ref{qrbar}
we have
\begin{equation}
\int\limits_{S_{N}^4}
\exp\Bigl(
N\alpha(R_{1,2}-R_{3,4})^2+\sum_{l\leq 4}(\vsi^l,\vec{v})
\Bigr)d\lambda_N^4
=
\int\limits_{S_{N}^4}
\exp\Bigl(
N\alpha(R_{1,2}-R_{3,4})^2+\sum_{l\leq 4}\eps_l |\vec{v}|
\Bigr)d\lambda_N^4.
\label{above}
\end{equation}
Since
\begin{eqnarray*}
(R_{1,2}-R_{3,4})^2 
&\leq& 2(\hR_{1,2}-\hR_{3,4})^2+ 
\frac{2}{N^2}(\eps_1\eps_2-\eps_3\eps_4)^2
\\
&\leq&
4\hR_{1,2}^2+4\hR_{3,4}^2+
\frac{2}{N^2}(\eps_1\eps_2-\eps_3\eps_4)^2,
\end{eqnarray*}
using (\ref{Fub2}), the right hand side of (\ref{above}) is bounded by
$$
a_N^4\int\limits_{[-\sqrt{N},\sqrt{N}]^4}
\exp\Bigl(
\frac{2\alpha}{N}(\eps_1\eps_2-\eps_3\eps_4)^2+\sum_{l\leq 4}\eps_l |\vec{v}|
\Bigr)d\vec{\eps}
\Bigl(
\int\limits_{S_{N-1}^2}\exp\bigl(
4\alpha N\hR_{1,2}^2
\bigr)d\lambda_{N-1}^2(\hS^1,\hS^2)
\Bigr)^2,
$$
where $d\vec{\eps}=d\eps_1\ldots d\eps_4.$
For a fixed $\hS^2\in S_{N-1},$ let $Q$ be an orthogonal transformation
in $\Reals^{N-1}$ such that 
$$
Q\hS^2=(0,\ldots,0,|\hS^2|)=(0,\ldots,0,\sqrt{N-1}).
$$
Then 
$$
\hR_{1,2}=\frac{1}{N-1}\Bigl(Q\hS^1,Q\hS^2\Bigr)
=\frac{1}{\sqrt{N-1}}(Q\hS^1)_{N-1}.
$$
Therefore, by rotational invariance and then (\ref{Fub2}), 
\begin{eqnarray*}
&&
\int\limits_{S_{N-1}^2}\exp\bigl(
4\alpha N\hR_{1,2}^2
\bigr)d\lambda_{N-1}^2(\hS^1,\hS^2)
=
\int\limits_{S_{N-1}}
\exp4\alpha\frac{N}{N-1}\eps^2 d\lambda_{N-1}(\hS)
\\
&&
\leq
a_{N-1}\int\limits_{-\sqrt{N-1}}^{\sqrt{N-1}}
\exp (5\alpha \eps^2)
\Bigl(
1-\frac{\eps^2}{N-1}
\Bigr)^{\frac{N-4}{2}}d\eps
\leq
L\int\limits_{-\infty}^{\infty}
\exp\bigl(5\alpha\eps^2 - L\eps^2\bigr)d\eps
\leq L
\end{eqnarray*}
for small enough $\alpha.$ Therefore, the right hand side
of (\ref{above}) is bounded for small $\alpha$ by
$$
L \int\limits_{[-\sqrt{N},\sqrt{N}]^4}
\exp\Bigl(
\frac{2\alpha}{N}(\eps_1\eps_2-\eps_3\eps_4)^2+\sum_{l\leq 4}\eps_l |\vec{v}|
\Bigr)d\vec{\eps}.
$$
Making the change of variables $\eps_l=\sqrt{N}x_l$ (as in (\ref{upthere})) proves that
$\e\bla \exp N \alpha(R_{1,2}-\bq)^{2}\bra_0$ is bounded up to a constant by
\begin{eqnarray*} 
&&
\e\int\limits_{[-1,1]^4} 
\exp N \Phi(\vec{x})d\vec{x}
\Bigr/
\Bigl(\int\limits_{-1}^{1} \exp N\varphi(x) dx\Bigr)^{4}
\end{eqnarray*}
where
\begin{equation}
\Phi(\vec{x})
=\Phi(x_1,x_2,x_3,x_4)=2\alpha(x_1 x_2 - x_3 x_4)^2+\sum_{l\leq 4}\varphi(x_l)
\label{Phi}
\end{equation}
and where $\varphi(x)$ was defined in (\ref{fi}).
We will use this bound only on the event $\{c\leq c'\}$
since by (\ref{cc})
$$
\e\bla \exp N \alpha(R_{1,2}-\bq)^{2}\bra_0 
\leq \exp(4N\alpha-LN)+
\e\bla \exp N \alpha(R_{1,2}-\bq)^{2}\bra_0 I(c\leq c')
$$
and $L$ can be made as large as necessary by taking $c'$ sufficiently large.
Since by (\ref{derr}), $\varphi''(x)\leq -1,$ for small enough $\alpha$
the function $\Phi(\vec{x})$ will be strictly concave on $[-1,1]^4.$ It is
obvious that for $\vec{x}_0=(x_0,x_0,x_0,x_0)$
$$
\frac{\partial \Phi}{\partial x_l}(\vec{x}_0)=\varphi'(x_0)=0
$$
which implies that $\vec{x}_0$ is the unique maximum of $\Phi.$
Strict concavity now implies
$$
\Phi(\vec{x})\leq 4\varphi(x_0)-\frac{1}{L}\sum_{l\leq 4}(x_l-x_0)^2
$$
and, thus,
\begin{eqnarray*}
&&
\int\limits_{[-1,1]^4}\exp N\Phi(\vec{x})d\vec{x}
\leq
\exp 4N\varphi(x_0)
\Bigl(
\int\limits_{-1}^{1}\exp\Bigl(-\frac{1}{L}N(x-x_0)^2\Bigr)dx
\Bigr)^4
\leq
\frac{L}{N^2}\exp 4N\varphi(x_0).
\end{eqnarray*}
Combining this with (\ref{Zlow}) finally proves that
$
\e\bla \exp N \alpha(R_{1,2}-\bq)^{2}\bra_0 \leq L.
$
The proof of the corresponding statement for 
$R_1-\bar{r}$ is similar.
\qed

\end{document}